# Construction of stationary trajectories for a model of a system of $N$ particles with interaction

Igor Pavlov

Abstract. For the classical $N-body$ problem, an approach is proposed based on the introduction of some natural in the physical sense optimization problems of mathematical programming for finding a conditional minimum for the characteristics of the system on the set of its possible states. The solution of these problems then makes it possible to construct families of flat stationary and periodic trajectories of the system and also to find relationships and estimates for the characteristics of the system on these trajectories. It is shown that when the system moves on a plane on trajectories generated by the global minimum in these optimization problems, at any time the minimum possible size of the system is achieved at each current level of its "cohesion" (or potential energy). Similar optimization problems are considered for finding a conditional minimum for the characteristics of a system in three-dimensional space. It is shown that the solution of these problems can be achieved only on flat trajectories of the system and is achieved, in particular, on the constructed flat stationary and periodic trajectories. In addition, it is shown that the trajectory of the system in three-dimensional space, at least at one point of which the minimum possible size of the system is achieved at the current value of its cohesion (or potential energy), can only be flat. And such trajectories are, in particular, flat stationary and periodic trajectories generated by the global minimum in the considered optimization problems.

## 1. Introduction

Consider a system of $N$ different material points (particles). Let $r_i$ be the radius vector and $m_i$ be the mass of the particle with index $i = 1, \ldots, N$. Further in sections $2-9$ by the state of the system at the moment of time $t \geq 0$ we will mean the set of radius vectors of all particles $r = (r_1, r_2, \ldots, r_N) \in \Re$, where $\Re$ is the set of all possible states (coordinate space):

$$\Re = \{r : |r_j - r_i| > 0, \ i \neq j\} \tag{1.1}$$

The strict inequalities in (1.1) correspond to the natural physical condition that the radius vectors of different particles do not coincide. (Here and below, vectors are denoted in bold, $|a|$ is the length of the vector $a$. The indices $i, j$ everywhere take the values $1, \ldots, N$.)

The interaction force $F_{ij} = F_{ij}(r)$, acting on the particle $i$ from the particle $j$ in the state $r = (r_1, r_2, \ldots, r_N) \in \Re$ is equal to

$$F_{ij}(r) = \frac{\gamma\, m_i\, m_j\, L_{ij}}{|L_{ij}|^3}, \quad i \neq j \tag{1.2}$$

Where $\gamma > 0$ is some constant, $L_{ij} = L_{ij}(r) = r_j - r_i$. The system is assumed to be closed, that is there are no external influences on it. The interaction force $F_i = F_i(r)$ acting on the particle $i$ from all other particles of the system in the state $r \in \Re$ has the form

$$F_i(r) = \sum_{j \neq i} F_{ij}(r) = \sum_{j \neq i} \frac{\gamma\, m_i\, m_j\, L_{ij}}{|L_{ij}|^3} \tag{1.3}$$

Formulas (1.2), (1.3) correspond to Newton's law of gravity, when the force of attraction between two particles with masses $m_i, m_j$ is inversely proportional to the square of the distance between them. This model corresponds to the known classical $N - body\ problem$.

The trajectory of the system $\boldsymbol{r}(t) = [\,\boldsymbol{r}_1(t),\ \boldsymbol{r}_2(t), ...,\ \boldsymbol{r}_N(t)] \in \Re,\ t \geq 0$ is given by a system of second-order differential equations

$$m_i\,\ddot{\boldsymbol{r}}_i = \boldsymbol{F}_i(\boldsymbol{r}_1,\ \boldsymbol{r}_2, ...,\ \boldsymbol{r}_N) \tag{1.4}$$

(where $\boldsymbol{r}_i = \boldsymbol{r}_i(t),\ \dot{\boldsymbol{r}}_i = \dot{\boldsymbol{r}}_i(t),\ \ddot{\boldsymbol{r}}_i = \ddot{\boldsymbol{r}}_i(t),\ i = 1,...,N$) with certain initial conditions

$$\boldsymbol{r}_i(0) = \boldsymbol{r}_{i0},\ \dot{\boldsymbol{r}}_i(0) = \boldsymbol{v}_{i0},\ i = 1, ..., N \tag{1.5}$$

A trajectory $\boldsymbol{r}(t) = [\,\boldsymbol{r}_1(t),\ \boldsymbol{r}_2(t), ...,\ \boldsymbol{r}_N(t)] \in \Re,\ t \geq 0$ satisfying the system of differential equations (1.4) under certain initial conditions (1.5) will be called stationary (or finite in terms of book [1]), if there exist constants $0 < C_1 \leq C_2 < \infty$ such that

$$C_1 \leq |\boldsymbol{r}_j(t) - \boldsymbol{r}_i(t)| \leq C_2 \quad \text{for each}\ \ i \neq j,\ \ t \geq 0,$$

In accordance with this definition, the stationarity (finiteness) of the trajectory of the system $\boldsymbol{r}(t) \in \Re,\ t \geq 0$ means that on this trajectory the system (and, possibly, its structure) is preserved in time in some compact form for all $t \geq 0$.

Below, in sections $2 - 9$, this problem is considered first on the plane. An approach is proposed based on the introduction of some natural in the physical sense optimization problems of mathematical programming for finding a conditional minimum for the characteristics of size (of compactness) and the "cohesion" of the system on the set of its possible states $\Re$. The solution of these problems then makes it possible to construct families of flat stationary and periodic trajectories of the system and also to find relationships and estimates for the characteristics of the system on these trajectories.

In sections 4 and 8, it is shown that when the system moves on a plane on trajectories generated by the global minimum in these optimization problems, at any time $t \geq 0$ the minimum possible size of the system is achieved at each current level of its "cohesion" (or potential energy). In section 5, lower bounds are obtained for the main characteristics of the system on these trajectories. Section 7 shows that local minima in these optimization problems also generate families of flat stationary and periodic trajectories of the system with similar relations and estimates for its characteristics. In section 9, the question of the structure of the system on the constructed flat trajectories is considered.

In section 10, similar optimization problems are considered for finding a conditional minimum for the characteristics of a system in three-dimensional space. It is shown (Theorems $10.1 - 10.3$) that the solution of these problems is achieved (for all $t \geq 0$) on the constructed flat stationary and periodic trajectories.

Section 11 shows (Theorem 11.1) that the trajectory of the system in three-dimensional space, at least at one point of which the minimum possible size of the system is achieved at the current value of its cohesion (or potential energy), can only be flat. And such trajectories are, in particular, flat stationary and periodic trajectories generated by the global minimum in the considered optimization problems.

## 2. Problems for finding a conditional minimum for system characteristics on a plane.

Further in sections $2 - 9$ this problem is first considered on a plane, that is, it is assumed that all particles of the system are in the same plane, and the position of the particles of the system on this plane is given by two-dimensional radius vectors $r_i = (x_i, y_i)$, $i = 1, ..., N$. The state $r$ of the system at time $t \geq 0$ is given by the set $r = (r_1, r_2, ..., r_N) \in \Re \subset E_{2N}$ of the radius vectors of all particles (where $E_n$ denotes the $n$-dimensional Euclidean space). Consider the state function (force function)

$$f(r) = \frac{1}{2} \sum_i \sum_{j \neq i} \frac{\gamma \, m_i \, m_j}{|r_j - r_i|} \,, \quad r \in \Re \tag{2.1}$$

coinciding up to a sign and a constant with the potential energy of the system $\Pi(r)$ in this model: $f(r) = -\Pi(r) + C$, [1] − [3].

Function (2.1) has the meaning of the work or energy required to completely destroy a given system, that is, to transfer all its particles from a given state $r \in \Re$ (at zero initial particle velocities) to an infinite distance from each other. Therefore, it is natural to call the function $f(r)$ the cohesion characteristic or simply the "cohesion" of the system in the state $r = (r_1, r_2, ..., r_N) \in \Re$. The cohesion of the system $f(r)$, as it is seen from (2.1), generally speaking, should increase with decreasing size (increasing compactness) of the system.

Further, we will assume that the center of mass of the system coincides with the origin of coordinates, that is, $\sum_i m_i r_i = 0$ (further, this condition will also automatically follow from the solution of the problem). As a characteristic of the size of the system in the state $r \in \Re$, we take its "root-mean-square size" (or simply the "size") $b = b(r) = \sqrt{\sum_i \beta_i |r_i|^2}$, where $\beta_i = m_i / \sum_j m_j$, $i = 1, ..., N$.

To simplify the formulas below, as a characteristic of the size of the system in the state $r \in \Re$, we will also use the function

$$g = g(r) = \sum_i m_i |r_i|^2 \tag{2.2}$$

which is related to $b = b(r)$ by the relation $g = b^2 \sum_j m_j$. (Where the function $g = g(r)$ also has the meaning of the moment of inertia of the system in the state $r \in \Re$.)

Consider the following problem of finding a conditional minimum on the set of system states $\Re \subset E_{2N}$: it is required to find

$$min \; g(r), \quad r \in \Re \tag{2.3}$$

under the condition $f(r) = f$. That is, it is required to find the state of the system $r \in \Re$ with minimum size for a given fixed level of its cohesion $f$. Or in the dual formulation ([4] - [6], etc.) of this problem: it is required to find

$$h(g) = min \; f(r), \quad r \in \Re \tag{2.4}$$

under the condition $g(r) = g$. That is, it is required to find the state of the system $r \in \Re$ with minimum cohesion for a fixed its size.

We will also consider this problem in the unconditional formulation: it is required to find

$$min \; [f(r) + \lambda g(r)], \quad r \in \Re \tag{2.5}$$

where $\lambda > 0$ is the Lagrange multiplier.

It is easy to show further that for any $\lambda > 0$ the minimum (2.5) is reached at some point (state) $r = (r_1, r_2, ..., r_N) \in \Re$. Indeed, it follows directly from the definition of the functions $f(r)$, $g(r)$ that, for any fixed $\lambda > 0$, the objective function $f(r) + \lambda g(r)$ in (2.5) is continuous in a closed bounded domain of the form $H_C = \{r \in \Re : f(r) + \lambda g(r) \leq C\} \subset \Re$ where $C > 0$ is any sufficiently large constant such that $C > inf[f(r) + \lambda g(r)] \geq 0$. Whence, taking into account the definition of the set $H_C$, it follows that the minimum in (2.5) is certainly attained at some (perhaps not the only) point $r \in H_C \subset \Re$.

Any state in which the minimum in (2.5) is reached for a given $\lambda > 0$ will be denoted by $r_\lambda = (r_{1\lambda}, r_{2\lambda}, \ldots, r_{N\lambda}) \in \Re$. The necessary conditions for the minimum (2.5) at the point $r_\lambda \in \Re$ have the form

$$\frac{\partial f}{\partial r_i}(r_\lambda) = -\lambda \frac{\partial g}{\partial r_i}(r_\lambda), \quad i = 1, \ldots, N \tag{2.6}$$

For the left and right sides of these equations, the following relations hold

$$\frac{\partial f}{\partial r_i}(r_\lambda) = F_i(r_\lambda), \quad \frac{\partial g}{\partial r_i}(r_\lambda) = 2m_i \, r_{i\lambda}, \quad i = 1, \ldots, N$$

Thus, the necessary conditions for a minimum (2.5) in the state $r_\lambda = (r_{1\lambda}, r_{2\lambda}, \ldots, r_{N\lambda}) \in \Re$ can also be written in the form

$$F_i(r_\lambda) = -2\lambda m_i \, r_{i\lambda} \quad \text{for each} \quad i = 1, \ldots, N \tag{2.7}$$

where $F_i(r_\lambda)$ is the force acting on the particle with index $i$ in this state. Whence, summing the left and right sides of (2.6) over $i = 1, \ldots, N$ we find $\sum_i m_i \, r_{i\lambda} = 0$, that is, in the state $r_\lambda \in \Re$, in which the minimum of the objective function in (2.5) is reached, the center of mass of the system coincides with the origin.

Note that equations (2.6), (2.7) must hold not only in the state of the global minimum $r_\lambda \in \Re$ in (2.5), but also in the state of any local minimum. By virtue of (2.7), the structure of the system in any state in which the global or local minimum of the objective function in problem (2.5) is reached, forms a "central configuration" ([7], [8], etc.). The number of such configurations, therefore, is no less than the number of different local minima (up to equivalence relations) in optimization problems of the form (2.3) – (2.5) (see also section 7 below).

The solution of these optimization problems (and similar more general problems in sections 10 and 11) makes it possible to construct families of flat stationary and periodic trajectories of the system, as well as to find relations and estimates for the characteristics of the system on these trajectories. In addition, this makes it possible to establish the existence and type of trajectories on which the minimum possible size (maximum compactness) of the system is achieved when it moves on a plane or in three-dimensional space (sections $3 - 11$).

### 3. A family of flat stationary trajectories of the system.

Let $r_\lambda = (r_{1\lambda}, r_{2\lambda}, \ldots, r_{N\lambda}) \in \Re \subset E_{2N}$ be any state of the system in which the minimum of the objective function in (2.5) is reached, where $r_{i\lambda} = (x_{i\lambda}, y_{i\lambda})$, $x_{i\lambda} = |r_{i\lambda}| \cos \varphi_{i\lambda}$, $y_{i\lambda} = |r_{i\lambda}| \sin \varphi_{i\lambda}$, $i = 1, \ldots, N$. Let us further define the trajectory of the system motion in the coordinate space based on the state $r_\lambda \in \Re$ as follows:

$$r_\lambda(t) = (r_{1\lambda}(t), \ldots, r_{N\lambda}(t)) \in \Re \tag{3.1}$$

Where $r_{i\lambda}(t) = (x_{i\lambda}(t), y_{i\lambda}(t))$, $x_{i\lambda}(t) = |r_{i\lambda}| \cos(\varphi_{i\lambda} + \omega t)$, $y_{i\lambda}(t) = |r_{i\lambda}| \sin(\varphi_{i\lambda} + \omega t)$, $i = 1, \ldots, N$; where $t \geq 0$, $\omega = \sqrt{2\lambda}$. Thus, the state $r_\lambda \in \Re$, in which the minimum in (2.5) is reached, is the initial state of this trajectory, $r_\lambda = r_\lambda(0)$. The proof of the following theorems follows directly from the definition trajectories $r_\lambda(t)$.

**Theorem 3.1** *Any trajectory of the type $r_\lambda(t)$ of motion of the system is a solution to the system of differential equations (1.4) with the initial conditions $r_i(0) = r_{i\lambda}$, $\dot{r}_i(0) = (-\omega \, y_{i\lambda}, \omega \, x_{i\lambda})$, $i = 1, \ldots, N$.*

**Theorem 3.2** *Any trajectory of the type $r_\lambda(t)$ of motion of the system is stationary (in the sense of the definition of section 1) and periodic, and the following equations of the form (2.7), (2.6) hold at each of its points:*

$$\frac{\partial f}{\partial r_i}[r_\lambda(t)] = -\lambda \frac{\partial g}{\partial r_i}[r_\lambda(t)], \quad t \geq 0, \quad i = 1, \ldots, N$$

$$F_i[r_\lambda(t)] = -2\lambda m_i\, r_{i\lambda}(t), \quad t \geq 0, \quad i = 1, \ldots, N$$

Thus, on each trajectory $r_\lambda(t)$ the internal forces of interaction between the particles of the system are such that all its particles move on this trajectory in a "consistent way", namely, they rotate around the origin of coordinates (the center of mass of the system) with the same angular velocity $\omega = \sqrt{2\lambda}$.

On each trajectory $r_\lambda(t)$ the system has a constant structure, which is preserved in time and is determined by the initial state $r_\lambda = (r_{1\lambda}, r_{2\lambda}, \ldots, r_{N\lambda})$ in which the minimum (2.5) is reached. And this state essentially depends on the number of particles $N$ in the system and the set of particle masses $m_1, m_2, \ldots, m_N$. The set of all such trajectories for different $\lambda > 0$ forms a one-parameter family of flat stationary and periodic trajectories of motion of the system:

$$r_\lambda(t) \in \mathfrak{R} \subset E_{2N}, \quad t \geq 0; \quad \lambda > 0 \qquad (3.2)$$

The existence (for any $N > 1$, $m_1, m_2, \ldots, m_N$) of such flat stationary and periodic trajectories is a consequence of the fact that the minimum of the objective function in (2.5) is achieved on the coordinate space $\mathfrak{R} \subset E_{2N}$ (see section 2) and the necessary conditions (2.7) for this minimum.

On each trajectory of the form $r_\lambda(t)$ the characteristics $f(r_\lambda(t))$ and $g(r_\lambda(t))$ of the system are constant and for them we will use the abbreviated notations $f_\lambda = f(r_\lambda(t))$, $g_\lambda = g(r_\lambda(t))$, etc.

### 4. Relations for system characteristics on flat stationary trajectories

**Theorem 4.1** *The solution of the optimization problem (2.4) — the function $h(g)$ is monotonically decreasing and differentiable with respect to the parameter of this problem $g > 0$. In addition, the following relation for the differentials is satisfied*

$$dh(g) = -[h(g)/2](dg/g) \qquad (4.1)$$

*Proof.* It follows from the definition of the functions $f(r)$, $g(r)$ that these functions satisfy the relations $f(\propto r) = f(r)/\propto$, $g(\propto r) = \propto^2 g(r)$ for any $\propto > 0$, $r = (r_1, r_2, \ldots, r_N) \in \mathfrak{R}$, where $\propto r = (\propto r_1, \propto r_2, \ldots, \propto r_N) \in \mathfrak{R}$. Denote by $B(g) = \{r \in \mathfrak{R}: g(r) = g\} \subset \mathfrak{R}$ the subset of states $r \in \mathfrak{R}$, in which the value of the function $g(r)$ is fixed and equal to the constant $g > 0$. Let $r^{(0)} = \left(r_1^{(0)}, \ldots, r_N^{(0)}\right) \in B(g_0)$ be the state, in which the minimum in (2.4) is reached for a given fixed value of $g = g_0 > 0$, that is $h(g_0) = f(r_0)$.

Let us show that then, for any other $g = g_1 > 0$ the minimum in (2.4) under the condition $g(r) = g_1$ is reached in the state $r^{(1)} = \left(r_1^{(1)}, \ldots, r_N^{(1)}\right) = \propto r^{(0)} \in B(g_1)$, where $\propto = \sqrt{g_1/g_0}$. Indeed, let $u = (u_1, \ldots, u_N)$ be an arbitrary state in $B(g_1)$. Then the state $r = (u/\propto) \in B(g_0)$, since $g(r) = g(u/\propto) = g(u)/\propto^2 = g_1/\propto^2 = g_0$ and hence

$$f(u) = f(\propto r) = f(r)/\propto \geq f(r^{(0)})/\propto = f(\propto r^{(0)}) = f(r^{(1)})$$

i.e. minimum in (2.4) on the set $B(g_1)$ is reached in the state $\boldsymbol{r}^{(1)} = \propto \boldsymbol{r}^{(0)} \in B(g_1)$. Thus $h(g_1) = f(\boldsymbol{r}^{(1)})$. Let us increment $g$ to $g + \Delta g$ in (2.4) and set $g = g_0$, $g + \Delta g = g_1$. Then

$$h(g + \Delta g) = h(g_1) = f(\boldsymbol{r}^{(1)}) = f(\propto \boldsymbol{r}^{(0)}) = f(\boldsymbol{r}^{(0)})/\propto = h(g_0)/\propto = = h(g)/\propto$$

where $\propto = (g_1/g_0)^{1/2} = [(g + \Delta g)/g]^{1/2}$, whence after simple transformations we find

$$\Delta h = h(g + \Delta g) - h(g) = h(g)\left\{\left(1 + \frac{\Delta g}{g}\right)^{-1/2} - 1\right\}$$

Whence

$$\Delta h = -\left(\frac{h(g)}{2g}\right)\Delta g + o(\Delta g), \quad if \quad \Delta g \to 0 \qquad (4.2)$$

Consequently, the function $h(g)$ decreases monotonically and is differentiable with respect to $g > 0$. From (4.2) relation (4.1) for the differentials follows. ∎

From (4.1), after integration, we find that the solution $h(g)$ of the problem of finding a conditional minimum in (2.4), depending on the parameter of this problem $g$, has the form

$$h(g) = C/\sqrt{g}, \qquad (4.3)$$

where the integration constant $C = C(\boldsymbol{m}) > 0$ is some function of the set of particle masses of the system $\boldsymbol{m} = (m_1, m_2, \dots, m_N)$.

Consider the set $D_{\boldsymbol{m}}$ of all possible values of the vector $(g, f)$ of characteristics of the system $g = g(\boldsymbol{r})$, $f = f(\boldsymbol{r})$ for a given set $\boldsymbol{m} = (m_1, m_2, \dots, m_N)$ of the masses of its particles:

$$D_{\boldsymbol{m}} = \{(g, f): g = g(\boldsymbol{r}), f = f(\boldsymbol{r}), \quad \boldsymbol{r} \in \Re\} \qquad (4.4)$$

From (4.3), (4.4) and the definition of the function $h(g)$ as a solution to the optimization problem (2.4), it follows that this set has the form

$$D_{\boldsymbol{m}} = \{(g, f): \quad f \geq C(\boldsymbol{m})/\sqrt{g}\} \qquad (4.5)$$

and its lower bound has the form

$$G_{\boldsymbol{m}} = \{(g, f): \quad f\sqrt{g} = C(\boldsymbol{m})\} \qquad (4.6)$$

**Theorem 4.2** *On each flat stationary and periodic trajectory $\boldsymbol{r}_\lambda(t)$ the "cohesion" of the system $f_\lambda = f(\boldsymbol{r}_\lambda(t))$ and the characteristic of its size $g_\lambda = g(\boldsymbol{r}_\lambda(t))$ are constant and are related by the relation*

$$f_\lambda = C(\boldsymbol{m})/\sqrt{g_\lambda} \qquad (4.7)$$

*Proof.* Let $H(\lambda)$ be the solution of the above problem of finding the minimum in (2.5) for a given $\lambda > 0$, i.e.

$$H(\lambda) = min \; [f(\boldsymbol{r}) + \lambda g(\boldsymbol{r})]$$

where the minimum is taken over all $\boldsymbol{r} \in \Re$. Then, taking into account the definition of the set $D_{\boldsymbol{m}}$ in (4.4), the function $H(\lambda)$ can also be represented in the form

$$H(\lambda) = min \; (f + \lambda g) \qquad (4.8)$$

where the minimum is taken over all $(g, f) \in D_{\boldsymbol{m}}$. It is easy to see that for an $\lambda > 0$ this minimum is attained at the lower boundary $G_{\boldsymbol{m}}$ of the region $D_{\boldsymbol{m}}$ given in (4.6). Calculating the minimum in (4.8) under the condition $(g,f) \in G_m \subset D_m$, after simple transformations we find, that the point $(g_\lambda, f_\lambda) \in G_m$, at which this minimum is reached for a given $\lambda > 0$, is determined by the equalities

$$g_\lambda = \left(C(\boldsymbol{m})/2\lambda\right)^{2/3}, \quad f_\lambda = [C(\boldsymbol{m})]^{2/3} (2\lambda)^{1/3} \quad (4.9)$$

Whence follows relation (4.7) between the characteristics $f_\lambda$, $g_\lambda$ of the system on flat stationary trajectories of the form $\boldsymbol{r}_\lambda(t)$. ∎

It follows from Theorem 4.2 that on each trajectory $\boldsymbol{r}_\lambda(t)$ the vector of characteristics of the system $(g, f)$ is always located on the lower boundary $G_m$ of the region $D_m$ of all possible values of this vector for a given system with a set of particle masses $\boldsymbol{m} = (m_1, m_2, \dots, m_N)$.

Taking into account the definition of dual optimization problems (2.3), (2.4), the minimum possible value of the characteristic $g$ of the system at a given level of its cohesion $f$ is determined from the equation $h(g) = f$. From the definition of problem (2.4), formula (4.3) and Theorem 4.2, taking into account also that $g = b^2 \sum_j m_j$, it follows

**Theorem 4.3** *When the system moves on a plane, the minimum possible root-mean-square size of the system $b$ at a given level of its cohesion $f$ is equal to*

$$b = \frac{C(\boldsymbol{m})}{f \sqrt{\sum_j m_j}} \quad (4.10)$$

*and this minimum possible size (and the maximum compactness in this sense) of the system is achieved on each trajectory $\boldsymbol{r}_\lambda(t)$ from the family (3.2).*

Further in section 8 it is shown (Theorem 8.3) that a similar fact holds for a more general two-parameter family of flat stationary trajectories generated by the global minimum $\boldsymbol{r}_\lambda \in \mathfrak{R} \subset E_{2N}$ in optimization problems of the form (2.3) – (2.5).

Section 11 also shows (Theorem 11.1) that the minimum root-mean-square size (and maximum compactness in this sense) of a system when it moves in three-dimensional space can be achieved only on flat trajectories, and such trajectories are, in particular, flat stationary and periodic trajectories $\boldsymbol{r}_\lambda(t)$ from the family (3.2) and some others flat trajectories.

On each trajectory $\boldsymbol{r}_\lambda(t) = (\boldsymbol{r}_{1\lambda}(t), \dots, \boldsymbol{r}_{N\lambda}(t))$ in accordance with its definition, the kinetic energy of the system $T$ is constant and

$$T = \sum_i m_i |\boldsymbol{v}_{i\lambda}(t)|^2 / 2 \equiv \omega^2 \sum_i m_i |\boldsymbol{r}_{i\lambda}|^2 / 2 = \omega^2 g_\lambda / 2 = \lambda g_\lambda \quad (4.11)$$

where $\boldsymbol{v}_{i\lambda}(t) = \dot{\boldsymbol{r}}_{i\lambda}(t)$, $\omega = \sqrt{2\lambda}$. Let $L = |\boldsymbol{L}|$ be the modulus of the angular momentum of the system $\boldsymbol{L}$ (constant due to known conservation laws). On the trajectory $\boldsymbol{r}_\lambda(t)$ lying in the plane $z = 0$, the angular momentum of the system $\boldsymbol{L} = (0, 0, L)$, where by definition of this trajectory

$$L = \sum_i m_i |\boldsymbol{v}_{i\lambda}(t)| |\boldsymbol{r}_{i\lambda}(t)| \equiv \omega \sum_i m_i |\boldsymbol{r}_{i\lambda}|^2 = \omega g_\lambda = \sqrt{2\lambda} g_\lambda \quad (4.12)$$

Then from (4.3), (4.9), (4.11), (4.12), taking into account also that $g = b^2 \sum_j m_j$, we find that on each trajectory $\boldsymbol{r}_\lambda(t)$ the characteristics of the system $f, T, L, \omega$ are constant and are related to its root-mean-square size $b$ (which is also constant) by the relations

$$f = \frac{C(\boldsymbol{m})}{b \sqrt{\sum_j m_j}}, \quad T = \frac{C(\boldsymbol{m})}{2b \sqrt{\sum_j m_j}}, \quad \omega = \frac{\sqrt{C(\boldsymbol{m})}}{b^{3/2} \left(\sum_j m_j\right)^{3/4}} \quad (4.13)$$

$$L = \sqrt{b}\,\sqrt{C(\boldsymbol{m})}\left(\sum_j m_j\right)^{1/4} \tag{4.14}$$

It also follows from formulas (4.13), (4.14) that on each trajectory $\boldsymbol{r}_\lambda(t)$ the characteristics of the system satisfy the relations

$$TL^2 = C^2(\boldsymbol{m})/2, \quad T\sqrt{g} = C(\boldsymbol{m})/2, \quad \omega^2 g^{3/2} = C(\boldsymbol{m}) \tag{4.15}$$

where $C(\boldsymbol{m})$ is some constant for a system with a given set of masses $\boldsymbol{m} = (m_1, m_2, \ldots, m_N)$ of particles. In addition, relations (4.13) imply that on each trajectory $\boldsymbol{r}_\lambda(t)$ the equality for the characteristics of the system is satisfied:

$$f = 2T \tag{4.16}$$

which also follows from the known virial theorem ([9], Sec. 5.5, etc.), according to which $\ddot{g} = 4T - 2f$, taking into account that on the trajectories $\boldsymbol{r}_\lambda(t)$, as follows from their definition, $g \equiv const$, $\ddot{g} \equiv 0$.

Further in section 7 it is shown that Theorems 4.2, 4.4 and relations (4.13) − (4.16) are also valid for all flat stationary trajectories of a system of the form (3.1), generated by any local minimum $r_\lambda^{(k)} \in \Re \subset E_{2N}$ in problems of the form (2.3) – (2.5), up to the replacement of the global minimum $\boldsymbol{r}_\lambda = (\boldsymbol{r}_{1\lambda}, \boldsymbol{r}_{2\lambda}, \ldots, \boldsymbol{r}_{N\lambda})$ by the local minimum $\boldsymbol{r}_\lambda^{(k)} = (\boldsymbol{r}_{1\lambda}^{(k)}, \boldsymbol{r}_{2\lambda}^{(k)}, \ldots, \boldsymbol{r}_{N\lambda}^{(k)})$ and the function $C(\boldsymbol{m})$ on the $C_k(\boldsymbol{m})$, where $C_k(\boldsymbol{m}) \geq C(\boldsymbol{m})$, $k = 1, \ldots, n$; $n$ is the number of local minima (see Section 7 below).

For the function $C(\boldsymbol{m})$, which is rather complicated in the general case, a lower estimate is obtained in the next Section 5.

## 5. Lower bounds for system characteristics on flat stationary trajectories

In accordance with Section 2, the minimum in problem (2.5) can be achieved only on a subset of states $G_0 = \{\boldsymbol{r} \in \Re: \sum_j m_j \boldsymbol{r}_j = 0\} \subset \Re \subset E_{2N}$. Whence it follows that the solution of the problem of finding the conditional minimum in (2.3), (2.4) is also achieved on this subset. Let us show that on the subset $G_0 \subset \Re$ the system characteristic $g(\boldsymbol{r}) = \sum_i m_i |\boldsymbol{r}_i|^2$ coincides with the function $g_1(\boldsymbol{r}) = B^{-1} \sum_{i<j} m_i m_j |\boldsymbol{L}_{ij}|^2$ where $B = \sum_i m_i$ is the mass of all particles of the system, $\boldsymbol{L}_{ij} = \boldsymbol{L}_{ij}(\boldsymbol{r}) = \boldsymbol{r}_j - \boldsymbol{r}_i$.

The partial derivatives of the function $g_1(\boldsymbol{r})$ for any $\boldsymbol{r} \in \Re$ have the form

$$\frac{\partial g_1}{\partial \boldsymbol{r}_k} = B^{-1} \frac{\partial}{\partial \boldsymbol{r}_k} \sum_{i<j} m_i m_j |\boldsymbol{L}_{ij}|^2 = B^{-1} m_k \sum_{j \neq k} m_j \frac{\partial}{\partial \boldsymbol{r}_k} |\boldsymbol{L}_{kj}|^2$$

Whence

$$\frac{\partial g_1}{\partial \boldsymbol{r}_k} = -2B^{-1} m_k \sum_{j \neq k} m_j \boldsymbol{L}_{kj}, \quad k = 1, \ldots, N; \quad \boldsymbol{r} \in \Re$$

where for the right side of this formula, the relations are fulfilled

$$\sum_{j \neq k} m_j \boldsymbol{L}_{kj} = \sum_{j \neq k} m_j (\boldsymbol{r}_j - \boldsymbol{r}_k) = \sum_j m_j \boldsymbol{r}_j - m_k \boldsymbol{r}_k - \boldsymbol{r}_k \sum_{j \neq k} m_j = -B\boldsymbol{r}_k$$

for each $\boldsymbol{r} \in G_0 \subset \Re$. Thus

$$\frac{\partial g_1}{\partial \boldsymbol{r}_k}(\boldsymbol{r}) = \frac{\partial g}{\partial \boldsymbol{r}_k}(\boldsymbol{r}) = 2m_k \boldsymbol{r}_k \quad \text{for each} \quad k = 1, \ldots, N; \quad \boldsymbol{r} \in G_0$$

whence follows the equality of the total differentials $dg_1(r) = dg(r)$, $r \in G_0$ and therefore the functions $g(r)$, $g_1(r)$ coincide everywhere on the set $G_0 \subset \Re$ up to a constant, that is $g(r) = g_1(r) + C$ for each $r \in G_0 \subset \Re$. Whence at $r_j \to 0$, $j = 1, \ldots, N$ it follows that $C = 0$ and $g(r) = g_1(r)$ for each $r \in G_0 \subset \Re$.

Thus, the problem of finding the conditional minimum (2.4) can further be represented in the following form: it is required to find

$$h(g) = \min \sum_{i<j} \frac{\gamma\, m_i\, m_j}{L_{ij}} \qquad (5.1)$$

where the minimum is taken under the conditions

$$B^{-1} \sum_{i<j} m_i m_j L_{ij}^2 = g \qquad (5.2)$$

$$\sum_j m_j x_j = 0, \quad \sum_j m_j y_j = 0 \qquad (5.3)$$

where the variables $L_{ij} = |L_{ij}|$ must satisfy additional constraints

$$L_{ij}^2 = (x_j - x_i)^2 + (y_j - y_i)^2, \quad i < j, \; (i,j) \in (1, \ldots, N) \qquad (5.4)$$

where $(x_i, y_i)$, $i = 1, \ldots, N$ form a set of $2N$ independent variables.

In accordance with Theorem 4.1 and formula (4.3) (see Section 4 above), the solution of this problem has the form $h(g) = C(\boldsymbol{m})/\sqrt{g}$, where $C(\boldsymbol{m})$ is some function of masses $\boldsymbol{m} = (m_1, m_2, \ldots, m_N)$ of particles of the system. Let us further consider a simpler problem: It is required to find the minimum of the same objective function

$$h_\mathrm{L}(g) = \min \sum_{i<j} \frac{\gamma\, m_i\, m_j}{L_{ij}} \qquad (5.5)$$

where the minimum is taken over all $L_{ij} > 0$, $i < j$, $(i,j) \in (1, \ldots, N)$ under the condition

$$B^{-1} \sum_{i<j} m_i m_j L_{ij}^2 = g \qquad (5.6)$$

but without additional restrictions (5.3), (5.4). Then the solution of this problem — function $h_\mathrm{L}(g)$ gives a lower estimate for the exact solution $h(g) = C(\boldsymbol{m})/\sqrt{g}$. Calculating the minimum in (5.5) under condition (5.6), after simple transformations we find $h_\mathrm{L}(g) = C_\mathrm{L}(\boldsymbol{m})/\sqrt{g}$,

$$C_\mathrm{L}(\boldsymbol{m}) = \gamma \left( \sum_{i<j} m_i m_j \right)^{3/2} \left( \sum_i m_i \right)^{-1/2} \qquad (5.7)$$

Whence, taking into account that $h_\mathrm{L}(g) \leq h(g)$, it follows

$$C_\mathrm{L}(\boldsymbol{m}) \leq C(\boldsymbol{m}) \qquad (5.8)$$

Thus $C_\mathrm{L}(\boldsymbol{m})$ gives a lower bound for the function $C(\boldsymbol{m}) = C(m_1, m_2, \ldots, m_N)$ (for any $N > 1$, $m_1, m_2, \ldots, m_N$). In the case $N = 3$, the lower bound (5.7) coincides with the exact formula for $C(\boldsymbol{m})$ (see Section 6 below).

From the lower estimate (5.7), taking into account the expressions (4.13), (4.14) obtained above, we then find the following estimates for the characteristics $f, T, L, \omega$ of the system on flat stationary and periodic trajectories $r_\lambda(t)$:

$$T \geq \frac{\gamma}{2b} \left( \sum_{i<j} m_i m_j \right)^{3/2} \left( \sum_i m_i \right)^{-1} \quad (5.9)$$

$$L \geq \sqrt{\gamma b} \left( \sum_{i<j} m_i m_j \right)^{3/4} \quad (5.10)$$

$$\omega \geq \frac{\sqrt{\gamma}}{b^{3/2}} \left( \sum_{i<j} m_i m_j \right)^{3/4} \left( \sum_i m_i \right)^{-1} \quad (5.11)$$

where $b$ is the root-mean-square size of the system. (Note that these estimates, as well as the relations for the system characteristics in (4.13), obtained in this model in the framework of classical mechanics and assuming zero particle sizes, cease to be valid if $b \to 0$)

### 6. Analytical solution for the case $N = 3$

Let us consider for the case $N = 3$ the solution of the problem (on the plane) for finding a minimum in the form (2.5): it is required to find

$$\min [f(r) + \lambda g(r)], \quad r \in \mathfrak{R} \subset E_{2N} \quad (6.1)$$

Further, it will be seen that in this case the solution of this problem, that is the state $r_\lambda = (r_{1\lambda}, r_{2\lambda}, r_{3\lambda})$ in which the minimum (6.1) is reached, generates a flat stationary and periodic trajectory of the system, coinciding with the known solution of Lagrange for the three-body problem ([3] Section 14, [9] Section 5.8).

For any vector $\boldsymbol{\alpha}$ we denote below by $\alpha$ the length of this vector, $\alpha = |\boldsymbol{\alpha}|$. Let $r_i = (x_i, y_i)$ be the radius-vector of the particle $i$, where $x_i, y_i$ are its Cartesian coordinates, and let $r_i, \varphi_i$ be its polar coordinates, where $r_i = |r_i|$, $x_i = r_i \cos\varphi_i$, $y_i = r_i \sin\varphi_i$, $i = 1,2,3$. Without loss of generality, we can assume that

$$0 = \varphi_1 < \varphi_2 < \varphi_3 < 2\pi \quad (6.2)$$

In the state $r = (r_1, r_2, r_3)$, in which the minimum in (6.1) is reached, the center of mass of the system coincides with the origin of coordinates (see item 2 above):

$$m_1 x_1 + m_2 x_2 + m_3 x_3 = 0 \quad (6.3)$$

$$m_1 y_1 + m_2 y_2 + m_3 y_3 = 0 \quad (6.4)$$

Since, by virtue of (6.2), $\varphi_1 = 0$, then

$$y_1 = 0, \quad x_1 = r_1 > 0 \quad (6.5)$$

Denote by $\mathbf{p}_{23} = (\alpha, \beta)$ the center of mass of particles with indices 2, 3:

$$\alpha = \frac{m_2 x_2 + m_3 x_3}{m_2 + m_3}, \quad \beta = \frac{m_2 y_2 + m_3 y_3}{m_2 + m_3}$$

By virtue of (6.3) – (6.5) $\beta = 0$, and

$$\alpha = -\frac{m_1 x_1}{m_2 + m_3} = -\frac{m_1 r_1}{m_2 + m_3} \quad (6.6)$$

Denote $L_{ij} = |\boldsymbol{L}_{ij}|$ the length of the vector $\boldsymbol{L}_{ij} = \boldsymbol{r}_j - \boldsymbol{r}_i$, $i \neq j$. From (6.5) follows
$$L_{12}^2 = (y_2 - y_1)^2 + (x_2 - x_1)^2 = y_2^2 + (x_2 - r_1)^2 =$$
$$= r_2^2 - x_2^2 + (x_2 - r_1)^2 = r_2^2 + r_1^2 - 2r_1 x_2 \tag{6.7}$$

Let us introduce the vectors $\boldsymbol{b}_2 = \boldsymbol{p}_{23} - \boldsymbol{r}_2$, $\boldsymbol{b}_3 = \boldsymbol{r}_3 - \boldsymbol{p}_{23}$. Taking into account the definition of the vector $\boldsymbol{p}_{23}$, the relation $m_2 \boldsymbol{b}_2 = m_3 \boldsymbol{b}_3$ is satisfied and

$$m_2 b_2 = m_3 b_3 \tag{6.8}$$

Whence follows the equality

$$L_{23}^2 = (b_2 + b_3)^2 = b_2^2 \left(\frac{m_2 + m_3}{m_3}\right)^2 \tag{6.9}$$

The value $b_2 = |\boldsymbol{b}_2|$ satisfies the relation

$$b_2^2 = |\boldsymbol{p}_{23} - \boldsymbol{r}_2|^2 = (\alpha - x_2)^2 + (\beta - y_2)^2 = (\alpha - x_2)^2 + y_2^2 \tag{6.10}$$

Whence, taking into account (6.6), we find

$$b_2^2 = \left(\frac{m_1 r_1}{m_2 + m_3} + x_2\right)^2 + r_2^2 - x_2^2 = r_2^2 + \left(\frac{m_1 r_1}{m_2 + m_3}\right)^2 + \frac{2 m_1 r_1 x_2}{m_2 + m_3} \tag{6.11}$$

From (6.9) − (6.11) follows

$$L_{23}^2 = r_2^2 \left(\frac{m_2 + m_3}{m_3}\right)^2 + r_1^2 \left(\frac{m_1}{m_3}\right)^2 + 2 r_1 x_2 \left(\frac{m_1 m_2 + m_1 m_3}{m_3^2}\right) \tag{6.12}$$

From (6.3) − (6.5), (6.8), taking into account that $|y_3/y_2| = b_3/b_2 = m_2/m_3$, after simple transformations follows

$$L_{31}^2 = (y_3 - y_1)^2 + (x_3 - x_1)^2 = y_2^2 \left(\frac{m_2}{m_3}\right)^2 + \left(r_1 \frac{m_1 + m_3}{m_3} + x_2 \frac{m_2}{m_3}\right)^2$$

Whence, taking into account that $y_2^2 = r_2^2 - x_2^2$, it follows

$$L_{31}^2 = r_2^2 \left(\frac{m_2}{m_3}\right)^2 + r_1^2 \left(\frac{m_1 + m_3}{m_3}\right)^2 + 2 r_1 x_2 \left(\frac{m_1 m_2 + m_2 m_3}{m_3^2}\right) \tag{6.13}$$

On the other hand, from these relations we find

$$r_3^2 = x_3^2 + y_3^2 = \left(x_1 \frac{m_1}{m_3} + x_2 \frac{m_2}{m_3}\right)^2 + y_2^2 \left(\frac{m_2}{m_3}\right)^2 =$$

$$= r_1^2 \left(\frac{m_1}{m_3}\right)^2 + r_2^2 \left(\frac{m_2}{m_3}\right)^2 + 2 r_1 x_2 \left(\frac{m_1 m_2}{m_3^2}\right)$$

Whence follows

$$2 r_1 x_2 = -r_1^2 \left(\frac{m_1}{m_2}\right) - r_2^2 \left(\frac{m_2}{m_1}\right) + r_3^2 \left(\frac{m_3^2}{m_1 m_2}\right) \tag{6.14}$$

From (6.7), (6.12) − (6.14) we further find the relations that connect the variables $r_i = |\boldsymbol{r}_i|$ and $L_{ij} = |\boldsymbol{L}_{ij}|$

$$L_{12}^2 = r_1^2 \left(\frac{m_1 + m_2}{m_2}\right) + r_2^2 \left(\frac{m_1 + m_2}{m_1}\right) - r_3^2 \left(\frac{m_3^2}{m_1 m_2}\right) \tag{6.15}$$

$$L_{23}^2 = -r_1{}^2\left(\frac{m_1{}^2}{m_2 m_3}\right) + r_2{}^2\left(\frac{m_2 + m_3}{m_3}\right) + r_3{}^2\left(\frac{m_2 + m_3}{m_2}\right) \tag{6.16}$$

$$L_{31}^2 = r_1{}^2\left(\frac{m_1 + m_3}{m_3}\right) - r_2{}^2\left(\frac{m_2{}^2}{m_1 m_3}\right) + r_3{}^2\left(\frac{m_1 + m_3}{m_1}\right) \tag{6.17}$$

Relations (6.15) − (6.17) are a system of linear equations with the determinant $D = (m_1 + m_2 + m_3)^3/(m_1 m_2 m_3) > 0$ with respect to the variables $r_1{}^2, r_2{}^2, r_3{}^2$. Whence it is not difficult to find also the system of inverse equations

$$\begin{aligned}r_1{}^2 &= L_{12}^2\left(\frac{m_2{}^2 + m_2 m_3}{(m_1 + m_2 + m_3)^2}\right) - L_{23}^2\left(\frac{m_2 m_3}{(m_1 + m_2 + m_3)^2}\right) \\ &+ L_{31}^2\left(\frac{m_2 m_3 + m_3{}^2}{(m_1 + m_2 + m_3)^2}\right)\end{aligned} \tag{6.18}$$

$$\begin{aligned}r_2{}^2 &= L_{12}^2\left(\frac{m_3 m_1 + m_1{}^2}{(m_1 + m_2 + m_3)^2}\right) + L_{23}^2\left(\frac{m_3{}^2 + m_3 m_1}{(m_1 + m_2 + m_3)^2}\right) \\ &- L_{31}^2\left(\frac{m_3 m_1}{(m_1 + m_2 + m_3)^2}\right)\end{aligned} \tag{6.19}$$

$$\begin{aligned}r_3{}^2 &= -L_{12}^2\left(\frac{m_1 m_2}{(m_1 + m_2 + m_3)^2}\right) + L_{23}^2\left(\frac{m_1 m_2 + m_2{}^2}{(m_1 + m_2 + m_3)^2}\right) \\ &+ L_{31}^2\left(\frac{m_1{}^2 + m_1 m_2}{(m_1 + m_2 + m_3)^2}\right)\end{aligned} \tag{6.20}$$

The objective function in (6.1) in this case has the form

$$f(\mathbf{r}) + \lambda g(\mathbf{r}) = \frac{\gamma\, m_1 m_2}{L_{12}} + \frac{\gamma\, m_2 m_3}{L_{23}} + \frac{\gamma\, m_3 m_1}{L_{31}} + \lambda(m_1 r_1{}^2 + m_2 r_2{}^2 + m_3 r_3{}^2) \tag{6.21}$$

where $L_{12}, L_{23}, L_{31}$ are functions of the variables $r_1, r_2, r_3$, which are defined in (6.15) – (6.17). From (6.21) follows

$$\frac{\partial}{\partial r_1}\{f(\mathbf{r}) + \lambda g(\mathbf{r})\} = -\frac{\gamma\, m_1 m_2}{L_{12}{}^2}\frac{\partial L_{12}}{\partial r_1} - \frac{\gamma\, m_2 m_3}{L_{23}{}^2}\frac{\partial L_{23}}{\partial r_1} - \frac{\gamma\, m_3 m_1}{L_{31}{}^2}\frac{\partial L_{31}}{\partial r_1} + 2\lambda m_1 r_1$$

Whence, taking into account (6.15) – (6.17), we find that in the state $\mathbf{r}_\lambda = (\mathbf{r}_{1\lambda}, \mathbf{r}_{2\lambda}, \mathbf{r}_{3\lambda})$, in which the minimum (6.1) is reached, the following relation is satisfied

$$\frac{\partial}{\partial r_1}\{f(\mathbf{r}_\lambda) + \lambda\, g(\mathbf{r}_\lambda)\} = -\frac{\gamma\, m_1(m_1 + m_2) r_{1\lambda}}{L_{12}{}^3} +$$

$$+ \frac{\gamma\, m_1{}^2 r_{1\lambda}}{L_{23}{}^3} - \frac{\gamma\, m_1(m_1 + m_3) r_{1\lambda}}{L_{31}{}^3} + 2\lambda m_1 r_{1\lambda} = 0$$

where $r_{1\lambda} = |\mathbf{r}_{1\lambda}|$, $L_{ij} = |\mathbf{r}_{j\lambda} - \mathbf{r}_{i\lambda}|$. Similarly, we find that in the state $\mathbf{r}_\lambda = (\mathbf{r}_{1\lambda}, \mathbf{r}_{2\lambda}, \mathbf{r}_{3\lambda})$ the following relations are also satisfied

$$\frac{\partial}{\partial r_2}\{f(\mathbf{r}_\lambda) + \lambda\, g(\mathbf{r}_\lambda)\} = -\frac{\gamma\, m_2(m_2 + m_1) r_{2\lambda}}{L_{12}{}^3} -$$

$$- \frac{\gamma\, m_2(m_2 + m_3) r_{2\lambda}}{L_{23}{}^3} + \frac{\gamma\, m_2{}^2 r_{2\lambda}}{L_{31}{}^3} + 2\lambda m_2 r_{2\lambda} = 0,$$

$$\frac{\partial}{\partial r_3}\{f(\mathbf{r}_\lambda) + \lambda\, g(\mathbf{r}_\lambda)\} = \frac{\gamma\, m_3{}^2 r_{3\lambda}}{L_{12}{}^3} -$$

$$-\frac{\gamma\, m_3(m_3+m_2)r_{3\lambda}}{L_{23}{}^3} - \frac{\gamma\, m_3\,(m_3+m_1)r_{3\lambda}}{L_{31}{}^3} + 2\lambda m_3 r_{3\lambda} = 0$$

where $r_{2\lambda} = |\boldsymbol{r}_{2\lambda}|$, $r_{3\lambda} = |\boldsymbol{r}_{3\lambda}|$. Whence follows a system of linear equations for the values of $(1/L_{12}{}^3)$, $(1/L_{23}{}^3)$, $(1/L_{31}{}^3)$:

$$\frac{m_1+m_2}{L_{12}{}^3} - \frac{m_1}{L_{23}{}^3} + \frac{m_3+m_1}{L_{31}{}^3} = \frac{2\lambda}{\gamma}$$

$$\frac{m_1+m_2}{L_{12}{}^3} + \frac{m_2+m_3}{L_{23}{}^3} - \frac{m_2}{L_{31}{}^3} = \frac{2\lambda}{\gamma}$$

$$-\frac{m_3}{L_{12}{}^3} + \frac{m_2+m_3}{L_{23}{}^3} + \frac{m_3+m_1}{L_{31}{}^3} = \frac{2\lambda}{\gamma}$$

This system has a unique solution

$$\frac{1}{L_{12}{}^3} = \frac{1}{L_{23}{}^3} = \frac{1}{L_{31}{}^3} = \frac{2\lambda}{\gamma(m_1+m_2+m_3)} \tag{6.22}$$

Thus, in the state $\boldsymbol{r}_\lambda = (\boldsymbol{r}_{1\lambda}, \boldsymbol{r}_{2\lambda}, \boldsymbol{r}_{3\lambda})$, in which the minimum of the objective function in (6.1) is reached (hence, also on the trajectory $\boldsymbol{r}_\lambda(t)$ generated by the state $\boldsymbol{r}_\lambda$), the particles of the system form an equilateral triangle with the distances between the particles

$$|L_{12}| = |L_{23}| = |L_{31}| = (\gamma B/2\lambda)^{1/3}$$

where $B = (m_1+m_2+m_3)$. And in accordance with section 3 on each trajectory $\boldsymbol{r}_\lambda(t)$, the system rotates as a whole around of the origin of coordinates with the angular velocity $\omega = \sqrt{2\lambda}$.

Thus, in the case $N = 3$ the solution $\boldsymbol{r}_\lambda = (\boldsymbol{r}_{1\lambda}, \boldsymbol{r}_{2\lambda}, \boldsymbol{r}_{3\lambda})$ of the problem of finding a minimum of the form (2.5) or (6.1) generates a flat stationary and periodic trajectory of the system, coinciding with the known solution of Lagrange for the three-body problem ([3] sec. 14, [9] sec. 5.8).

From equations (6.18) − (6.20), we find further that in the state of the minimum $\boldsymbol{r}_\lambda = (\boldsymbol{r}_{1\lambda}, \boldsymbol{r}_{2\lambda}, \boldsymbol{r}_{3\lambda})$ (hence, also on the trajectory $\boldsymbol{r}_\lambda(t)$ generated by this state), the following relations are fulfilled:

$$|\boldsymbol{r}_{1\lambda}| = \left(\frac{\gamma}{2\lambda}\right)^{1/3} \frac{(m_2{}^2 + m_2 m_3 + m_3{}^2)^{1/2}}{(m_1+m_2+m_3)^{2/3}} \tag{6.23}$$

$$|\boldsymbol{r}_{2\lambda}| = \left(\frac{\gamma}{2\lambda}\right)^{1/3} \frac{(m_1{}^2 + m_1 m_3 + m_3{}^2)^{1/2}}{(m_1+m_2+m_3)^{2/3}} \tag{6.24}$$

$$|\boldsymbol{r}_{3\lambda}| = \left(\frac{\gamma}{2\lambda}\right)^{1/3} \frac{(m_1{}^2 + m_1 m_2 + m_2{}^2)^{1/2}}{(m_1+m_2+m_3)^{2/3}} \tag{6.25}$$

From (6.22) − (6.25) it is not difficult to find further that the function $C(\boldsymbol{m})$, connecting the characteristics of the system (see sections 4 and 5 above), in this case has the form

$$C(\boldsymbol{m}) = \frac{\gamma(m_1 m_2 + m_2 m_3 + m_1 m_3)^{3/2}}{(m_1+m_2+m_3)^{1/2}}$$

and coincides in this case with the lower bound (5.7) obtained above in section 5.

## 7. Flat stationary trajectories generated by local minima

Below in this section it is shown that any local minimum in problems of the form (2.3) – (2.5) (considered on the plane) also generates a certain family of flat stationary and periodic trajectories of the system.

Note that if in the problems (2.3) – (2.5) the local minimum is reached in some state $r \in \Re \subset E_{2N}$, then taking into account the form of the functions $f(r)$, $g(r)$, it is also reached in any state of the form

$$r(\beta) \in \Re \subset E_{2N}, \quad 0 \leq \beta \leq 2\pi \tag{7.1}$$

where $r(\beta)$ the state obtained from $r$ by turning through the angle $\beta$. In this sense, all states of the form (7.1) are equivalent. Below different local minima at points $r^{(k)} \in \Re$ are understood as points (states) with different equivalence subsets of the form (7.1).

The proof of the following Theorem 7.1 for local minima in problems (2.3) – (2.5) is similar to the proof of Theorem 4.1 and Theorems 4.2, 4.3 for the case of the global minimum $r_\lambda \in \Re \subset E_{2N}$. Let the $B(g) = \{r \in \Re : g(r) = g\} \subset \Re \subset E_{2N}$ be a subset of points (states) of the coordinate space $\Re$ with a fixed value of the function $g(r) = g$. Let us consider the problem in the form (2.4) for finding the conditional minimum of the function $f(r)$.

**Theorem 7.1** *Let in the conditional minimum problem (2.4) for some fixed value of the parameter $g = g_0 > 0$ there exist $n \geq 1$ local minima of the function $f(r)$ on the set $B(g_0)$ at points $r_k(g_0) \in B(g_0)$, $k = 1, \ldots, n$. Then*

*1) For any other value of the parameter $g > 0$ there are also exist $n$ local minima of the function $f(r)$ on the set $B(g)$ at points $r_k(g) \in B(g)$, $k = 1, \ldots, n$, where the point of the local minimum*

$$r_k(g) = \alpha(g) r_k(g_0), \quad k = 1, \ldots, n \tag{7.2}$$

*where $\alpha(g) = \sqrt{(g/g_0)}$.*

*2) For any $g > 0$ the value of the k-th local minimum, that is, the function $h_k(g) = f[r_k(g)]$ has the form*

$$h_k(g) = C_k(m)/\sqrt{g}, \quad k = 1, \ldots, n \tag{7.3}$$

*where $C_k(m)$ is some function of the masses $m = (m_1, m_2, \ldots, m_N)$ of the particles of the system, such that $C_L(m) \leq C(m) \leq C_k(m)$, $k = 1, \ldots, n$, where $C_L(m)$ is the lower estimate (5.7) for the function $C(m)$ obtained above in Sections 4 and 5.*

*3) In the problem (2.5) for each $\lambda > 0$ the local minimum is reached at the point $r_\lambda^{(k)} = r_k[g_k(\lambda)]$, where $g_k(\lambda) = [C_k(m)/2\lambda]^{2/3}$, $k = 1, \ldots, n$.*

*Proof.* For any $g > 0$ and any point $r = (r_1, r_2, \ldots, r_N) \in B(g)$ denote by

$$V_\varepsilon(r) = \{u \in B(g) : |r_i - u_i| < \varepsilon, \quad i = 1, \ldots, N\} \subset B(g)$$

the $\varepsilon$-neighborhood of the point $r$ on the set $B(g)$. Since, according to the condition of the theorem, the function $f(r)$ on the set $B(g_0)$ has the local minimum at the point $r_k(g_0)$, then for some $\varepsilon > 0$ the following inequality is satisfied

$$f(\boldsymbol{u}) \geq f[\boldsymbol{r}_k(g_0)] \quad \text{for each} \quad \boldsymbol{u} \in V_\varepsilon[\boldsymbol{r}_k(g_0)], \quad k = 1,\ldots,n$$

Using the abbreviation $\alpha = \alpha(g) = \sqrt{(g/g_0)}$, we determine the value of $\delta$ from the condition $\delta = \alpha\varepsilon$. Let $\boldsymbol{w} \in V_\delta[\boldsymbol{r}_k(g)]$ be any point from the $\delta$-neighborhood of $\boldsymbol{r}_k(g)$ on the set $B(g)$. Then the point $\boldsymbol{u} = (\boldsymbol{w}/\alpha) \in B(g_0)$, since $g(\boldsymbol{u}) = g(\boldsymbol{w}/\alpha) = g(\boldsymbol{w})/\alpha^2 = g/\alpha^2 = g_0$, In addition, due to the equality of $\delta = \alpha\varepsilon$ the point $\boldsymbol{u} = (\boldsymbol{w}/\alpha) \in V_\varepsilon[\boldsymbol{r}_k(g_0)]$. Whence, taking into account also (7.4), follows

$$f(\boldsymbol{w}) = f(\alpha \boldsymbol{u}) = f(\boldsymbol{u})/\alpha \geq f[\boldsymbol{r}_k(g_0)]/\alpha = f[\alpha \boldsymbol{r}_k(g_0)] = f[\boldsymbol{r}_k(g)]$$

for all $\boldsymbol{w} \in V_\delta[\boldsymbol{r}_k(g)]$. Consequently, at the point $\boldsymbol{r}_k(g) \in B(g)$ the function $f(\boldsymbol{r})$ has a local minimum on the set $B(g)$, $k = 1,\ldots,n$, which proves the first assertion of Theorem 7.1.

We put further $g_0 = g + \Delta g$. Then from the definition of the function $h_k(g) = f[\boldsymbol{r}_k(g)]$, taking into account also (7.2), it follows

$$\Delta h_k = h_k(g + \Delta g) - h_k(g) = f[\boldsymbol{r}_k(g_0)] - f[\boldsymbol{r}_k(g)] =$$
$$= f[\boldsymbol{r}_k(g)/\alpha] - f[\boldsymbol{r}_k(g)] = f[\boldsymbol{r}_k(g)](\alpha - 1) = h_k(g)(\alpha - 1)$$

where $\alpha = (g/g_0)^{1/2} = [(g + \Delta g)/g]^{-1/2}$, whence similarly to the proof of Theorem 4.1, the relation for the differentials follows

$$dh_k(g) = -[h_k(g)/2g]dg \tag{7.5}$$

Whence follows the expression (7.3) for the function $h_k(g)$, where $C_k = C_k(\boldsymbol{m})$ is the integration constant in (7.5), $k = 1,\ldots,n$. And any local minimum $h_k(g) \geq h(g)$, where $h(g) = C(\boldsymbol{m})/\sqrt{g}$ is the global minimum in this problem, obtained above in Theorem 4.1 and formula (4.3). Therefore $C_\mathrm{L}(\boldsymbol{m}) \leq C(\boldsymbol{m}) \leq C_k(\boldsymbol{m})$, $k = 1,\ldots,n$. Where $C_\mathrm{L}(\boldsymbol{m})$ is the lower estimate (5.7) for the function $C(\boldsymbol{m})$, which implies the second assertion of Theorem 7.1.

Thus, in the problem (2.4) the number of local minima $\boldsymbol{r}_k(g)$ on the set of $B(g) \subset \Re \subset E_{2N}$ is the same for all values of the problem parameter $g > 0$. At the point of the local minimum $\boldsymbol{r}_k(g) \in B(g)$ the Lagrange relations are satisfied

$$\frac{\partial f}{\partial \boldsymbol{r}_i} = -\lambda_k \frac{\partial g}{\partial \boldsymbol{r}_i}, \quad i = 1,\ldots,N \tag{7.6}$$

where $\lambda_k$ is the Lagrange multiplier. In addition, at this point the equality for the total differentials of these functions $df = -\lambda_k dg$ is satisfied, where $k = 1,\ldots,n$. Whence, taking into account (7.3) and the definition of the function $h_k(g) = f[\boldsymbol{r}_k(g)]$, follows

$$\lambda_k = -\frac{dh_k}{dg} = \frac{C_k(\boldsymbol{m})}{2g\sqrt{g}} \tag{7.7}$$

Whence, assuming $g_k(\lambda) = [C_k(\boldsymbol{m})/2\lambda]^{2/3}$, we obtain that the local minimum is attained at the point $\boldsymbol{r}_k[g_k(\lambda)]$ in the problem (2.5) for given $\lambda > 0$, $k = 1,\ldots,n$, which proves the third assertion of the theorem 7.1. ∎

Using the abbreviated notation $\boldsymbol{r}_\lambda^{(k)} = \boldsymbol{r}_k[g_k(\lambda)]$, from (7.6) and (7.7), we find that for each $\lambda > 0$ at the point $\boldsymbol{r}_\lambda^{(k)}$ the following relations are fulfilled

$$\frac{\partial f(\boldsymbol{r}_\lambda^{(k)})}{\partial \boldsymbol{r}_i} = -\lambda \frac{\partial g(\boldsymbol{r}_\lambda^{(k)})}{\partial \boldsymbol{r}_i}, \quad i = 1,\ldots,N; \quad k = 1,\ldots,n \tag{7.8}$$

Consequently, in each state of the local minimum $r_\lambda^{(k)} = (r_{1\lambda}^{(k)}, r_{2\lambda}^{(k)}, \ldots, r_{N\lambda}^{(k)})$ the following equations similar to (2.6) are satisfied:

$$F_i(r_\lambda^{(k)}) = -2\lambda m_i r_{i\lambda}^{(k)}, \qquad i = 1, \ldots, N$$

for all $\lambda > 0$, $k = 1, \ldots, n$. Whence it follows that each local minimum $r_\lambda^{(k)} \in \Re$ in problem (2.5) generates (similarly to Section 3) a flat stationary trajectory $r_\lambda^{(k)}(t) \in \Re$ of the system of the form (3.1) with the initial state of the trajectory – the local minimum $r_\lambda^{(k)} = (r_{1\lambda}^{(k)}, r_{2\lambda}^{(k)}, \ldots, r_{N\lambda}^{(k)}) \in \Re$ instead of the global minimum $r_\lambda = (r_{1\lambda}, r_{2\lambda}, \ldots, r_{N\lambda}) \in \Re$.

Consequently, the set of all local minima $r_\lambda^{(k)} \in \Re \subset E_{2N}$ for various $\lambda > 0$, $k = 1, \ldots, n$ generates $n$ of various families of flat stationary trajectories of the system

$$r_\lambda^{(k)}(t) = \left[r_{1\lambda}^{(k)}(t), r_{2\lambda}^{(k)}(t), \ldots, r_{N\lambda}^{(k)}(t)\right] \in \Re, \quad t \geq 0\,; \quad \lambda > 0 \qquad (7.9)$$

similar to the family (3.2) of trajectories $r_\lambda(t) \in \Re, t \geq 0$ generated by the global minimum $r_\lambda \in \Re$ in problem (2.5).

Proceeding from Theorem 7.1, it is easy to show, repeating the reasoning of Section 4, that the relations and estimates obtained above for the characteristics of the system on trajectories $r_\lambda(t)$ from the family (3.2) remain valid also for trajectories $r_\lambda^{(k)}(t)$ generated by each local minimum, replacing the function $C(m)$ to $C_k(m)$, $k = 1, \ldots, n$.

From Theorem 7.1, taking into account (7.2) and the definition of the state $r_\lambda^{(k)} = r_k[g_k(\lambda)]$, it follows that if in problem (2.5) the local minimum $r_{\lambda_0}^{(k)}$ exists for the value of the parameter $\lambda = \lambda_0 > 0$, then for any other $\lambda > 0$ the local minimum $r_\lambda^{(k)}$ also exists and $r_\lambda^{(k)} = (\lambda_0/\lambda)^{1/3} r_{\lambda_0}^{(k)}$ for each $k = 1, \ldots, n$. Consequently, for different values of the parameter $\lambda > 0$ the local minimum $r_\lambda^{(k)}$ in problem (2.5) is on a straight line ("ray") $r_\lambda^{(k)} = \alpha(\lambda) r_{\lambda_0}^{(k)}$ in the coordinate space $\Re \subset E_{2N}$, where $\alpha(\lambda) = (\lambda_0/\lambda)^{1/3}$.

Whence it follows that on any trajectory $r_\lambda^{(k)}(t)$ from the k-th family in (7.9) the state (and, consequently, the structure) of the system at any time $t \geq 0$ coincides with some fixed state $r_{\lambda_0}^{(k)}$, if we multiply the radius-vectors of all particles by the constant $\alpha(\lambda) = (\lambda_0/\lambda)^{1/3}$ and turn to angle $\beta = \omega t = \sqrt{2\lambda}\, t$. Consequently, all states of the system on trajectories from each family in (7.9) are equivalent up to rotation and scale.

Thus, local minima in optimization problems of the form (2.3) − (2.5) generate families (7.9) of flat stationary and periodic trajectories of the system, similar to the family (3.2) of trajectories $r_\lambda(t)$ for the case of a global minimum with similar relations (see above Sections 4 and 5) for the characteristics of the system on these trajectories.

But if a global minimum (perhaps not the only one) always exists, then the question of the existence and number of local minima in these optimization problems, depending on the dimension $N$ of the problem and the masses of the particles $m_1, m_2, \ldots, m_N$ of the system, remains unclear. Nevertheless, for a given set of parameters $N$, $m_1, m_2, \ldots, m_N$, $\lambda > 0$ one or another local minimum in problem (2.5) (and, consequently, the structure of the system on the trajectory generated by this minimum) can be found numerically by the gradient descent method (see also section 9 below).

## 8. A family of flat stationary trajectories of a system with pulsating size

In this section, we show that each (global or local) minimum in optimization problems of the form (2.3) − (2.5) (considered on the plane) generates a two-parameter family of flat stationary and periodic trajectories of the system with "pulsating size". Some relations and estimates for the characteristics of the system on these trajectories are obtained. It is shown (Theorem 8.4) that on all trajectories generated by the global minimum in problems of this type, at each point of the trajectory the minimum possible root-mean-square size (and the maximum compactness in this sense) of the system is achieved for each current level of its cohesion (or potential energy).

Let $r_\lambda = (r_{1\lambda}, r_{2\lambda}, \ldots, r_{N\lambda}) \in \Re \subset E_{2N}$ be the state in which the global minimum is reached in problem (2.5). Then, in accordance with section 2, this state satisfies the system of equations

$$F_j(r_\lambda) = -2\lambda m_j \, r_{j\lambda}, \quad j = 1, \ldots, N \tag{8.1}$$

or in expanded form

$$\sum_{k \neq j} \frac{\gamma \, m_j \, m_k (r_{k\lambda} - r_{j\lambda})}{|r_{k\lambda} - r_{j\lambda}|^3} = -2\lambda m_j \, r_{j\lambda}, \quad j = 1, \ldots, N \tag{8.2}$$

In state $r_\lambda = (r_{1\lambda}, r_{2\lambda}, \ldots, r_{N\lambda})$ each radius-vector $r_{j\lambda} = (x_{j\lambda}, y_{j\lambda})$ of the particle with index j let associate with a complex value $z_{j\lambda} = x_{j\lambda} + i y_{j\lambda}$, $j = 1, \ldots, N$. Then, by virtue (8.2), the set of $z_{j\lambda}$, $j = 1, \ldots, N$ satisfies the system of complex equations

$$\sum_{k \neq j} \frac{\gamma \, m_j \, m_k (z_{k\lambda} - z_{j\lambda})}{|z_{k\lambda} - z_{j\lambda}|^3} = -2\lambda m_j \, z_{j\lambda}, \quad j = 1, \ldots, N \tag{8.3}$$

For each particle with index j, we introduce a complex function

$$Z_{j\lambda}(t) = z_{j\lambda} a(t) e^{i\varphi(t)}, \quad j = 1, \ldots, N \tag{8.4}$$

where $a(t) > 0$, $\varphi(t)$ are real functions that do not depend on the index $j = 1, \ldots, N$, that is, these functions are the same for all particles of the system. What does it mean stretching (at $a(t) > 1$) or compression (at $a(t) < 1$) and rotation of the system by the angle $\varphi(t)$ at the time $t \geq 0$. Due to (8.4), the second derivative of the function $Z_{j\lambda} = Z_{j\lambda}(t)$ has the form

$$\ddot{Z}_{j\lambda} = z_{j\lambda} e^{i\varphi}(\ddot{a} - a\dot{\varphi}^2) + i \, z_{j\lambda} e^{i\varphi}(2\dot{a}\dot{\varphi} + a\ddot{\varphi}) \tag{8.5}$$

From the system of equations (8.2), (8.3), taking into account (8.4), (8.5), it is easy to find that the set of complex functions $Z_{j\lambda} = Z_{j\lambda}(t), j = 1, \ldots, N$ satisfies for all $t \geq 0$ the system of equations

$$\sum_{k \neq j} \frac{\gamma \, m_j \, m_k (Z_{k\lambda} - Z_{j\lambda})}{|Z_{k\lambda} - Z_{j\lambda}|^3} = -2\lambda m_j \, \ddot{Z}_{j\lambda}, \quad j = 1, \ldots, N \tag{8.6}$$

if the real functions $a = a(t)$, $\varphi = \varphi(t)$ satisfy the differential equations

$$a^2(\ddot{a} - a\dot{\varphi}^2) = -2\lambda \tag{8.7}$$

$$2\dot{a}\dot{\varphi} = -a\ddot{\varphi} \tag{8.8}$$

which are similar, up to a constant on the right-hand side of (8.7), to the well-known equations for the two-body problem ([1], etc.). Equation (8.8) implies

$$a^2(t)\dot{\varphi}(t) = C \tag{8.9}$$

which corresponds to Kepler's law of constant sectoral velocity known for the two-body problem ([1], etc.). From (8.7) - (8.9) after the transition to polar coordinates $\varphi$, $a = a(\varphi)$, further follows equation $(d_2 u/d\varphi^2) + u = (2\lambda/C^2)$, where $u = 1/a$, $a = a(\varphi)$. Whence, for a given fixed $\lambda > 0$, follows similarly to the well-known solution for the two-body problem [1]

$$a(\varphi) = \frac{p}{1 - e\cos\varphi} \qquad (8.10)$$

where $p = a_0(1 - e)$ and the eccentricity $e = 1 - (a_0^3 \omega_0^2 / 2\lambda)$. Where the parameters $a_0 > 0$, $\omega_0 \neq 0$ are determined by the initial conditions $a_0 = a(0)$, $\omega_0 = \dot\varphi(0)$ for $t = 0$ in equations (8.7), (8.8). We further assume that these parameters satisfy the inequalities $0 < a_0^3 \omega_0^2/(2\lambda) < 2$, that is $|e| < 1$. Without loss of generality, we further assume, that $\varphi(0) = 0$. Then, taking into account (8.10), the equality $\dot a(0) = 0$ is automatically satisfied. In accordance with (8.9), (8.10). The angular coordinate $\varphi = \varphi(t)$ is related to $t \geq 0$ by the expression

$$Ct = \int_0^\varphi \frac{p^2 d\alpha}{(1 - e\cos\alpha)^2} \qquad (8.11)$$

where taking into account (8.9) $C = a_0^2 \omega_0$. The parameter $\omega_0 = \dot\varphi(0)$ has the mean of the angular velocity of rotation of the system at the initial moment $t = 0$, that is, in the initial state $r_\lambda = (r_{1\lambda}, r_{2\lambda}, \ldots, r_{N\lambda}) \in \Re$ that generates the trajectory of the system of the form (8.4) or (8.12).

Thus, for a given fixed $\lambda > 0$ the state $r_\lambda = (r_{1\lambda}, r_{2\lambda}, \ldots, r_{N\lambda})$ of the minimum in (2.5) generates a two-parameter (with parameters $a_0$, $\omega_0$) family of flat trajectories of the system:

$$\tilde r_\lambda(t, a_0, \omega_0) = \{ \tilde r_{1\lambda}(t), \tilde r_{1\lambda}(t), \ldots, \tilde r_{N\lambda}(t) \}, \quad t \geq 0 \qquad (8.12)$$

where $\tilde r_{j\lambda}(t) = [\tilde x_{j\lambda}(t), \tilde y_{j\lambda}(t)]$,

$$\tilde x_{j\lambda}(t) = |r_{j\lambda}| \, a[\varphi(t)] \cos\varphi(t), \quad j = 1, \ldots, N \qquad (8.13)$$

$$\tilde y_{j\lambda}(t) = |r_{j\lambda}| \, a[\varphi(t)] \sin\varphi(t), \quad j = 1, \ldots, N \qquad (8.14)$$

Where $r_{j\lambda}$ is the radius vector of the particle with index j in the initial state $r_\lambda = (r_{1\lambda}, r_{2\lambda}, \ldots, r_{N\lambda})$ generating this family.

Any trajectory $\tilde r_\lambda(t, a_0, \omega_0)$ is stationary (in the sense of the definition in section 1) and periodic, if $|e| < 1$, that is the parameters $a_0, \omega_0$ satisfy the inequalities $0 < a_0^3 \omega_0^2/(2\lambda) < 2$. In the particular case, if $a_0 = 1$, $\omega_0 = \sqrt{2\lambda}$, then the trajectory $\tilde r_\lambda(t, a_0, \omega_0)$ coincides with the trajectory $r_\lambda(t)$ from the family (3.2) considered above.

On the trajectory $\tilde r_\lambda(t, a_0, \omega_0)$ each particle j of the system moves along an elliptic curve (in polar coordinates) of the form

$$|\tilde r_{j\lambda}(\varphi)| = |r_{j\lambda}| a(\varphi) = |r_{j\lambda}| \frac{a_0(1 - e)}{1 - e\cos\varphi}, \quad j = 1, \ldots, N \qquad (8.15)$$

And in this case the structure of the system remains constant and coincides with its structure in the initial state $r_\lambda = (r_{1\lambda}, r_{2\lambda}, \ldots, r_{N\lambda})$ up to the rotation of the system by the angle $\varphi = \varphi(t)$ and the scale, i. e. multiplying the length of the radius-vector of each particle $|r_{j\lambda}|$ by the value $a[\varphi(t)]$.

Let the $D = D(\lambda, a_0, \omega_0)$ be the class of all trajectories $\tilde r_\lambda(t, a_0, \omega_0)$ with parameters $\lambda > 0$, $a_0 > 0$, $\omega_0 \neq 0$ satisfying the inequalities

$$0 < a_0^3 \omega_0^2/(2\lambda) < 2. \qquad (8.16)$$

And let $D_1 = D(\lambda, 1, \omega_0) \subset D$ be a subclass of all such trajectories with parameters $\lambda > 0$, $a_0 = 1$, $\omega_0 \neq 0$ such that

$$\omega_0^2/(2\lambda) < 2 \tag{8.17}$$

Let us show that, in fact, the three-parameter class $D$ coincides with its two-parameter subclass $D_1$.

**Theorem 8.1** *The classes of trajectories $D$ and $D_1$ coincide, $D = D_1$.*

*Proof.* Let $r_{\lambda_0} \in \Re$ be the point at which the minimum (2.5) is reached for some fixed $\lambda_0$. The points $r_\lambda \in \Re$ and $r_{\lambda_0} \in \Re$ at which this minimum is reached for different values of the parameter $\lambda > 0$ are related by the relation $r_\lambda = (\lambda_0/\lambda)^{1/3} r_{\lambda_0}$ in accordance with Theorem 7.1 (see item 7 above). Let further some (arbitrary) trajectory in the class $D$ have parameters ($\lambda_0$, $a_0$, $\omega_0$), satisfying condition (8.15). Since $\varphi(0) = 0$, the initial state from which this trajectory starts has the form $a_0 r_{\lambda_0} \in \Re$. Then in the class $D_1$ there is a trajectory with parameters ($\lambda = \lambda_0/a_0^3$, $\omega_0$) that starts from the same state $r_\lambda = (\lambda_0/\lambda)^{1/3} r_{\lambda_0} = a_0 r_{\lambda_0} \in \Re$ and then coincides for all $t \geq 0$ with trajectory with parameters ($\lambda_0$, $a_0$, $\omega_0$) in the class $D$. Whence, taking into account also that $D_1 \subset D$, it follows that the classes $D$ and $D_1$ coincide. ∎

In accordance with this theorem, flat stationary and periodic trajectories from the class $D_1$ will be further denoted as $\tilde{r}_\lambda(t, \omega_0) = \tilde{r}_\lambda(t, 1, \omega_0)$, where the parameters $\lambda > 0$, $\omega_0 \neq 0$ of the trajectory satisfy inequality (8.17).

Let $\varphi = \varphi(t)$ be the angle of rotation of the system at the moment $t \geq 0$ on the trajectory $\tilde{r}_\lambda(t, \omega_0)$. It is easy to find the characteristics of the system $g = g(\varphi)$, $b = b(\varphi)$, $f = f(\varphi)$, $T = T(\varphi)$ for each $\varphi$ on this trajectory. In accordance with (8.15)

$$g(\varphi) = g_\lambda a^2(\varphi), \quad b(\varphi) = b_\lambda a(\varphi), \quad f(\varphi) = f_\lambda/a(\varphi) \tag{8.18}$$

where $g_\lambda = g(r_\lambda)$, $b_\lambda = b(r_\lambda)$, $f_\lambda = f(r_\lambda)$ are the values of these characteristics in the initial state $r_\lambda = (r_{1\lambda}, r_{2\lambda}, \ldots, r_{N\lambda})$ that generates this trajectory. Whence, taking into account (8.10), it follows that on each trajectory $\tilde{r}_\lambda(t, \omega_0)$ the root-mean-square size $b(\varphi)$ of the system periodically changes (pulsing) between the minimum and maximum values:

$$b_\lambda \frac{\omega_0^2}{2\lambda(1 + |e|)} \leq b(\varphi) \leq b_\lambda \frac{\omega_0^2}{2\lambda(1 - |e|)}$$

From (8.18) and equality (4.7) obtained above in Theorem 4.2, it follows

**Theorem 8.2** *On each trajectory $\tilde{r}_\lambda(t, \omega_0)$, the cohesion of the system $f = f(\varphi)$ and the characteristic $g = g(\varphi)$ satisfy the following relations*

$$f(\varphi)\sqrt{g(\varphi)} \equiv f_\lambda \sqrt{g_\lambda} = C(\boldsymbol{m}) \quad \text{for each } \varphi \geq 0 \tag{8.19}$$

which is a generalization of equality (4.7) for trajectories $r_\lambda(t)$ from the family (3.2). Where $C(\boldsymbol{m})$ is the function of masses $\boldsymbol{m} = (m_1, m_2, \ldots, m_N)$ of particles of the system introduced above in section 4.

By virtue of equality (8.19), the point $(g, f)$ of the characteristics of the system $g = g(\varphi)$, $f = f(\varphi)$ on each trajectory $\tilde{r}_\lambda(t, \omega_0)$ is always located on the lower boundary $G_{\boldsymbol{m}}$ of the set $D_{\boldsymbol{m}}$ in (4.4) of all possible values of the vector of characteristics $(g, f)$ for a given system with a set of particle masses $\boldsymbol{m} = (m_1, m_2, \ldots, m_N)$ (see section 4 above).

**Theorem 8.3** *On each trajectory $\tilde{r}_\lambda(t, \omega_0)$, the cohesion of the system $f = f(\varphi)$ and the kinetic energy of the system $T = T(\varphi)$ are satisfy to the following relations (for each $\varphi \geq 0$)*

$$f(\varphi) = \frac{C(m)}{b(\varphi)\sqrt{\sum_j m_j}} \tag{8.20}$$

$$T(\varphi) = \frac{C(m)}{2b(\varphi)\sqrt{\sum_j m_j}} \left(\frac{1 + e^2 - 2e\cos\varphi}{1 - e\cos\varphi}\right) \tag{8.21}$$

*where $b(\varphi)$ is the root-mean-square size of the system for a given value of $\varphi$.*

*Proof.* Equality (8.20) for the cohesion of the system $f(\varphi)$ follows directly from equality (8.19), also taking into account that $g(\varphi) = b^2(\varphi)\sum_j m_j$. In accordance with the known conservation laws ([1], sec. 6; [2], sec. 8) on each trajectory $\tilde{r}_\lambda(t, \omega_0)$ the total energy $E = T - f$ of the system is constant, i.e.

$$T(\varphi) - f(\varphi) = T(0) - f(0) \quad \text{for each } \varphi \geq 0 \tag{8.22}$$

Where the following relation is fulfilled

$$T(0) = \frac{1}{2}\sum_j m_j |z_{j\lambda}|^2 [\dot{a}(0)]^2 + \frac{1}{2}\sum_j m_j |z_{j\lambda}|^2 [a(0)]^2 [\dot{\varphi}(0)]^2$$

whence, учитывая (4.11), (4.16), (see section 4 above), as well as the equalities $\dot{a}(0) = 0$, $\dot{\varphi}(0) = \omega_0$, $a(0) = a_0 = 1$, we find

$$T(0) = \frac{\omega_0^2}{2}\sum_j m_j |r_{j\lambda}|^2 = \frac{\omega_0^2}{2} g(r_\lambda) = \lambda g_\lambda (1 - e) = f_\lambda (1 - e)/2$$

Where $e = 1 - (\omega_0^2/2\lambda)$. Whence, taking into account (8.10), (8.18), follows

$$T(\varphi) = T(0) + f(\varphi) - f(0) = \frac{f_\lambda(1 + e^2 - 2e\cos\varphi)}{2(1 - e)} \tag{8.23}$$

From (8.23), taking into account also (8.10), (8.18), (8.19), then follows (8.21). ∎

Using the lower estimate (5.7) for the function $C(m)$ obtained above in Section 5, from (8.21) we also find the lower estimate for the kinetic energy of the system $T(\varphi)$ on the any trajectory $\tilde{r}_\lambda(t, \omega_0)$ for each $\varphi \geq 0$

$$T(\varphi) \geq \frac{\gamma}{2b(\varphi)} \left(\sum_{i<j} m_i m_j\right)^{3/2} \left(\sum_i m_i\right)^{-1} \frac{1 + e^2 - 2e\cos\varphi}{1 - e\cos\varphi}$$

This estimate is a generalization of the lower estimate (5.9) for trajectories $r_\lambda(t)$ from the family (3.2) with a constant size of the system $b(\varphi) \equiv b_\lambda$.

From (8.18), (8.23) it is easy to find that the average values of the characteristics of the system $T(\varphi), f(\varphi)$ on the trajectory $\tilde{r}_\lambda(t, \omega_0)$ (along the angular coordinate $\varphi$) are equal

$$\langle T \rangle = \frac{1}{2\pi}\int_0^{2\pi} T(\varphi)\, d\varphi = \frac{f_\lambda(1 + e^2)}{2(1 - e)}$$

$$\langle f \rangle = \frac{1}{2\pi}\int_0^{2\pi} f(\varphi)\, d\varphi = f_\lambda/(1 - e) = \frac{2\langle T \rangle}{1 + e^2}$$

Whence, in particular, it follows that the maximum value of the average cohesion $\langle f \rangle = 2\langle T \rangle$ with respect to the average kinetic energy $\langle T \rangle$ takes place at the $e = 0$, that is on trajectories $r_\lambda(t)$ from the family (3.2) with a constant size $b(\varphi) \equiv b_\lambda$.

From (8.20), theorem 4.2 and formula (4.10) it follows

**Theorem 8.4**  *Let $\tilde{r}_\lambda(t, \omega_0)$ be any flat trajectory from the class $D_1$. Then, at each point of this trajectory (at any time $t \geq 0$), the minimum possible size of the system $b = b(\varphi)$ is achieved for each current level of the cohesion of the system $f = f(\varphi)$:*

$$b(\varphi) = \frac{C(m)}{f(\varphi)\sqrt{\sum_j m_j}} \qquad (8.24)$$

This theorem is actually a consequence of the fact that on each trajectory $\tilde{r}_\lambda(t, \omega_0)$, the vector of characteristics of the system $(g, f)$ is always located on the lower boundary $G_m$ of the set $D_m$ of all possible values $(g, f)$ for a given system with a set of particle masses $m = (m_1, m_2, \ldots, m_N)$ (see relation (8.19) and Section 4 above).

**Remark 1** on Theorem 8.4   The function $f = f(\varphi)$ coincides with the potential energy $\Pi = \Pi(\varphi)$ of the system up to a sign and an arbitrary constant, $f(\varphi) = -\Pi(\varphi) + C$. Thus, theorem 8.4 can also be formulated in terms of the potential energy of the system in the following equivalent form: *Let $\tilde{r}_\lambda(t, \omega_0)$ be any flat trajectory from the class $D_1$. Then, at each point of this trajectory, the minimum possible size of the system $b = b(\varphi)$ is achieved for each current level of the potential energy of the system $\Pi = \Pi(\varphi)$.*

**Remark 2** on Theorem 8.4   At each point of the trajectory $\tilde{r}_\lambda(t, \omega_0)$, the value $f(\varphi)$ is uniquely determined by the kinetic energy $T(\varphi)$ due to the equality $T(\varphi) - f(\varphi) \equiv E$, where the constant $E$ is the total energy of the system on this trajectory. Thus, Theorem 8.4 can also be formulated in the following equivalent form:   *Let $\tilde{r}_\lambda(t, \omega_0)$ be any flat trajectory from the class $D_1$. Then, at each point of this trajectory, the minimum possible size of the system $b = b(\varphi)$ is achieved for each current level of the kinetic energy of the system $T = T(\varphi)$.*

Thus, the state $r_\lambda = (r_{1\lambda}, r_{2\lambda}, \ldots, r_{N\lambda}) \in \Re \subset E_{2N}$, in which the global minimum in problem (2.5) is reached, generates a two-parameter family (with parameters $\lambda > 0$, $\omega_0 \neq 0$, $\omega_0^2/(2\lambda) < 2$) of flat stationary and periodic trajectories $\tilde{r}_\lambda(t, \omega_0)$ of the system with pulsing size $b(\varphi)$. (In the particular case, if the parameters $\omega_0 = \sqrt{2\lambda}$, this family coincides with the family (3.2) of trajectories $r_\lambda(t)$ with constant size $b(\varphi) \equiv b_\lambda$.)

Which corresponds (given the above region of parameter values  $\lambda > 0$, $\omega_0 \neq 0$, $\omega_0^2/(2\lambda) < 2$) to the following known general result of Arnold ([10],  p.130):   *"in n-body problem, there is a set of initial conditions that has a positive Lebesgue measure and such that if the initial positions and velocities of the bodies belong to this set, then the bodies always remain at a limited distance from each other. Astronomers have long expressed a similar hypothesis, but recent mathematicians who dealt with this issue, starting with Birkhoff, tended to the opposite opinion (see [11], [3])."*

It is essential that due to Theorem 8.4 on all trajectories $\tilde{r}_\lambda(t, \omega_0)$ generated by any global minimum in (2.5), the minimum possible size $b = b(\varphi)$ (and maximum compactness in this sense) of the system is achieved for each current level of its cohesion $f$ (or potential energy).

In the case $N = 3$ in the state $r_\lambda = (r_{1\lambda}, r_{2\lambda}, r_{3\lambda})$, in which the minimum  in (2.5) is reached, the particles of the system form an equilateral triangle (see Section 6 above), and in this case the family of trajectories  $\tilde{r}_\lambda(t, \omega_0)$ generated by the state $r_\lambda$, coincides with the known solution of Lagrange for the three-body problem ([3], sec. 14; [9], sec. 5.8).

Note that the above relations remain valid if  by $r_\lambda$ we mean the state of any local minimum $r_\lambda^{(k)} = (r_{1\lambda}^{(k)}, r_{2\lambda}^{(k)}, \ldots, r_{N\lambda}^{(k)}) \in \Re \subset E_{2N}$ in problem (2.5) up to replacing the global minimum $r_\lambda$ by the local minimum $r_\lambda^{(k)}$ and the function $C(m)$ with the function $C_k(m)$ (see Section 7 above). With the difference that on trajectories $\tilde{r}_\lambda(t, \omega_0)$ generated by one or another local minimum $r_\lambda^{(k)}$, the minimum possible size $b = b(\varphi)$ of the system in (8.24), generally speaking, is not achieved (if the local minimum of the objective function in (2.5) is strictly greater than the global minimum).

## 9. The structure of the system on flat stationary trajectories

The structure of the system at any moment in time $t \geq 0$ is determined by its state $r(t) = [r_1(t), r_2(t), \ldots, r_N(t)]$ at this moment. For a family (3.2) of flat stationary and periodic trajectories $r_\lambda(t)$, generated by the state $r_\lambda \in \Re \subset E_{2N}$, in which the minimum (2.5) is reached, the structure of the system at any time $t \geq 0$ remains constant. And this structure coincides with the structure of the initial state $r_\lambda = (r_{1\lambda}, r_{2\lambda}, \ldots, r_{N\lambda})$ up to rotation through the angle $\varphi = \omega t$, where $\omega = \sqrt{2\lambda}$.

For trajectories of a more general form $\tilde{r}_\lambda(t, \omega_0)$ (see above section 8) generated by the state of the global minimum $r_\lambda \in \Re \subset E_{2N}$, the structure of the system at any time also remains constant and coincides with the structure of the initial state $r_\lambda = (r_{1\lambda}, r_{2\lambda}, \ldots, r_{N\lambda})$ up to rotation and scale. Similarly, for trajectories of the form $r_\lambda^{(k)}(t)$, generated by the local minimum $r_\lambda^{(k)} \in \Re \subset E_{2N}$ in (2.5), the structure of the system is determined by the state $r_\lambda^{(k)} = (r_{1\lambda}^{(k)}, r_{2\lambda}^{(k)}, \ldots, r_{N\lambda}^{(k)})$, $k = 1, \ldots, n$ (see section 7 above).

In the case of $N = 3$, it follows from the previous formulas (6.22) – (6.25) that if $m_1 \geq m_2 \geq m_3$, then $|r_{1\lambda}| \leq |r_{2\lambda}| \leq |r_{3\lambda}|$. Thus, in the state $r_\lambda = (r_{1\lambda}, r_{2\lambda}, r_{3\lambda})$, in which the minimum (2.5) is reached (and on the trajectory $r_\lambda(t)$, generated by this state), particles with a smaller $m_i$ correspond to a greater distance $|r_{i\lambda}|$ from the origin of coordinates (the center of mass of the system).

Taking into account the general relations for the characteristics $f(r)$ and $g(r)$ of the system, it is natural to further assume that with an increase of the number $N$ of particles, a similar dependence holds. That is, in the state of a minimum $r_\lambda = (r_{1\lambda}, r_{2\lambda}, \ldots, r_{N\lambda}) \in \Re \subset E_{2N}$ in (2.5) and on the flat stationary trajectories generated by this minimum, particles with a smaller $m_i$ must also correspond to a greater distance $|r_{i\lambda}|$ from the origin of coordinates (the center of mass of the system). And such a dependence with an increase in the number of particles $N$ should lead to the structure of the system in the form of a spiral (with heavy core and light "tail") on these trajectories with the maximum compactness of the system (see also above Section 4, Section 8).

In the general case, for any $N$, the exact analytical solution of problems for finding a minimum in $(2.3) - (2.5)$ is rather complicated. Nevertheless, for given fixed parameters $N$, $m_1, m_2, \ldots, m_N$ these problems can be solved by numerical methods. The complexity of calculations in this case increases rapidly with increasing problem dimension. Examples obtained for relatively small $N$ using the gradient descent method ([4] - [6] etc.) for systems with a uniformly decreasing "spectrum" of particle masses $m_1, m_2, \ldots, m_N$ show the structure of the system in the form of a spiral on flat trajectories generated by these optimization problems.

## 10. Problems of finding a conditional minimum for the characteristics of a system in three-dimensional space.

In this section we consider some problems, natural in the physical sense, of finding a conditional minimum for the characteristics of a system in three-dimensional space, which are a generalization of the problems considered above for the characteristics of a system on a plane. Further we show that the solution of these problems is achieved (for each $t \geq 0$) on the flat stationary and periodic trajectories of the form $r_\lambda(t)$ and $r_\lambda^{(k)}(t)$ constructed above. The main theorem 10.1 shows that the solution of these problems of the first type can be achieved only on flat trajectories of the system. Proceeding from this theorem, in the next section it is shown (Theorem 11.1) that the minimum possible root-mean-square size (and the maximum compactness in this sense) of a system when it moves in three-dimensional space can be achieved only on flat trajectories.

Above, when constructing a family of flat stationary trajectories of the system, the state of the system was understood as a set of radius-vectors $r = (r_1, r_2, ..., r_N) \in \Re \subset E_{2N}$ of all particles of the system on the plane $(x, y)$, where $r_i = (x_i, y_i)$ is the radius vector of the particle with index $i = 1, ..., N$. Further, under the state of the system at the moment of time $t \geq 0$ in three-dimensional space $(x, y, z)$ we mean a pair $(r, v) \in \Omega = \Re \times \mho$, where $r = (r_1, r_2, ..., r_N)$ is the set of radius-vectors of the particles of the system in three-dimensional space, $r_i = (x_i, y_i, z_i)$ is the radius-vector of the particle with index $i$, $r \in \Re$, where

$$\Re = \{r : |r_j - r_i| > 0, \ i \neq j \} \subset E_{3N}, \tag{10.1}$$

$v = (v_1, v_2, ..., v_N)$ is a set of velocity vectors of various particles of the system, $v_i = (\alpha_i, \beta_i, \gamma_i)$ is the velocity vector of the particle with index $i$, $v \in \mho \subset E_{3N}$, where the set $\mho$ of possible values of $v = (v_1, v_2, ..., v_N)$ coincides with $E_{3N}$. The set (phase space) of all possible states $(r, v)$ of a system of $N$ particles in this model has view $\Omega = \Re \times \mho \subset E_{6N}$.

Let $L(r, v) = \sum_i [r_i, m_i v_i]$ be the angular momentum of the system (where $[r_i, m_i v_i]$ is the vector product of vectors $r_i, m_i v_i$), $T(v) = \sum_i m_i |v_i|^2 / 2$ is the kinetic energy of the system, and let $g(r) = \sum_i m_i |r_i|^2$ is the characteristic of the size of the system in the state $(r, v) \in \Omega$. The function $g(r)$ is related to the root-mean-square size of the system $b(r)$ by the same relation as in the problem on the plane: $g(r) = b^2(r) \sum_j m_j$.

Similar to the problems of section 2 for finding a minimum for the characteristics of a system on a plane, consider the following problem for finding a conditional minimum on the set of states $\Omega$ for system characteristics in three-dimensional space: it is required to find

$$\min g(r) \tag{10.2}$$

where the minimum is taken over all states $(r, v) \in \Omega$ under conditions on the characteristics of the system

$$L(r, v) = L \tag{10.3}$$

$$T(v) = T \tag{10.4}$$

where $L = (L_x, L_y, L_z)$, $T > 0$ are some constants. Condition (10.3) at the angular momentum of the system corresponds to three conditions on the components of this vector. We further assume without loss of generality that the vector $L$ is directed along the axis $z$, that is, $L = (0, 0, L)$ where $L \neq 0$.

In a more detailed form in terms of variables $(x_i, y_i, z_i)$, $(\alpha_i, \beta_i, \gamma_i)$, this problem has the following form: it is required to find

$$\min g(r) = \min \sum_i m_i (x_i^2 + y_i^2 + z_i^2) \tag{10.5}$$

where the minimum is taken over all states $(r, v) \in \Omega$ under the conditions

$$L_x(r, v) = \sum_i m_i (y_i \gamma_i - z_i \beta_i) = 0 \tag{10.6}$$

$$L_y(r, v) = \sum_i m_i (z_i \alpha_i - x_i \gamma_i) = 0 \tag{10.7}$$

$$L_z(\mathbf{r}, \mathbf{v}) = \sum_i m_i(x_i \beta_i - y_i \alpha_i) = L \qquad (10.8)$$

$$T(\mathbf{v}) = \frac{1}{2}\sum_i m_i(\alpha_i^2 + \beta_i^2 + \gamma_i^2) = T \qquad (10.9)$$

where $L \neq 0$, $T > 0$ are constants.

**Theorem 10.1** *The solution of the problem of finding a minimum in (10.2) under conditions (10.3), (10.4) can be achieved only on flat trajectories of motion of the system (that is, only at points of the phase space $(\mathbf{r}, \mathbf{v}) \in \Omega$ that belong to some flat trajectories).*

*Proof.* Let us consider the proof for the given problem presented in a more detailed form $(10.5) - (10.9)$. It suffices to show that in this problem the necessary conditions for a minimum at a point in the phase space $\mathbf{r}_i = (x_i, y_i, z_i)$, $\mathbf{v}_i = (\alpha_i, \beta_i, \gamma_i)$, $i = 1, \ldots, N$ are equalities $z_i = 0$, $\gamma_i = 0$ for each $i = 1, \ldots, N$.

The functions $g(\mathbf{r})$, $\mathbf{L}(\mathbf{r}, \mathbf{v})$, $T(\mathbf{v})$ are continuous, together with partial derivatives, on the closure of the phase space $\Omega^* \subset E_{6N}$, which in this case, taking into account the definition of the set $\mathfrak{R}$ in (10.1), coincides with the $E_{6N}$. Let $B \subset \Omega$ be a subset of the phase space over which the minimum in (10.5) is taken, that is, the set of all states $(\mathbf{r}, \mathbf{v}) \in \Omega$ for which conditions $(10.6) - (10.9)$ are satisfied. And let $B^* \subset \Omega^*$ be the closure of this set (that is, it also includes points that satisfy these conditions, but at which the radius-vectors of different particles $\mathbf{r}_i, \mathbf{r}_j$ can coincide).

Consider first the problem of finding the minimum (10.5) on the set $B^*$. It is easy to show that this minimum is reached at some (perhaps not the only) point $(\mathbf{r}, \mathbf{v}) \in B^*$. Indeed, let $B_C^* \subset B^*$ be a subset of all points of the phase space that belong to $B^*$ and such that $g(\mathbf{r}) \leq C$, where $C$ is any sufficiently large constant such that $C > \inf g(\mathbf{r})$ on the set $B^*$. The set $B_C^*$ is closed due to the continuity of the functions $\mathbf{L}(\mathbf{r}, \mathbf{v})$, $T(\mathbf{v})$, $g(\mathbf{r})$ and is bounded, taking into account also the condition (10.9) at the function $T(\mathbf{v})$. Thus, the function $g(\mathbf{r})$ is continuous on a closed bounded set $B_C^*$, whence it follows that the minimum in (10.5) is reached on the set $B_C^*$. Consequently, this minimum is also reached on the set $B^*$.

The necessary conditions for the minimum in (10.5) on the set $B^*$ at the point $(\mathbf{r}, \mathbf{v}) \in B^*$ have the form

$$\text{grad } g(\mathbf{r}) = a_1 \text{ grad } L_x(\mathbf{r}, \mathbf{v}) + a_2 \text{ grad } L_y(\mathbf{r}, \mathbf{v}) + a_3 \text{ grad } L_z(\mathbf{r}, \mathbf{v}) +$$
$$+ b \text{ grad } T(\mathbf{v}) \qquad (10.10)$$

where $a_1, a_2, a_3, b$ are the Lagrange multipliers. Or, taking into account that

$$\frac{\partial g}{\partial \mathbf{v}_i} = 0, \qquad \frac{\partial T}{\partial \mathbf{r}_i} = 0, \quad i = 1, \ldots, N,$$

equations (10.10) are further written as

$$\frac{\partial g}{\partial \mathbf{r}_i} = a_1 \frac{\partial L_x}{\partial \mathbf{r}_i} + a_2 \frac{\partial L_y}{\partial \mathbf{r}_i} + a_3 \frac{\partial L_z}{\partial \mathbf{r}_i}, \quad i = 1, \ldots, N \qquad (10.11)$$

$$0 = a_1 \frac{\partial L_x}{\partial \mathbf{v}_i} + a_2 \frac{\partial L_y}{\partial \mathbf{v}_i} + a_3 \frac{\partial L_z}{\partial \mathbf{v}_i} + b \frac{\partial T}{\partial \mathbf{v}_i}, \quad i = 1, \ldots, N \qquad (10.12)$$

Equations (10.11), (10.12) together with conditions $(10.6) - (10.9)$ form a system of $(6N + 4)$ equations for $(6N + 4)$ unknowns, namely, for the minimum point $\mathbf{r}_i = (x_i, y_i, z_i)$, $\boldsymbol{v}_i = (\alpha_i, \beta_i, \gamma_i)$, $i = 1,\ldots,N$; and coefficients $a_1, a_2, a_3, b$. Vector equations (10.11) in variables $(x_i, y_i, z_i)$, $(\alpha_i, \beta_i, \gamma_i)$, have the form

$$2x_i = -a_2 \gamma_i + a_3 \beta_i, \quad i = 1,\ldots,N \tag{10.13}$$

$$2y_i = a_1 \gamma_i - a_3 \alpha_i, \quad i = 1,\ldots,N \tag{10.14}$$

$$2z_i = -a_1 \beta_i + a_2 \alpha_i, \quad i = 1,\ldots,N \tag{10.15}$$

Similarly, equations (10.12) in variables $(x_i, y_i, z_i)$, $(\alpha_i, \beta_i, \gamma_i)$, are written as

$$0 = a_2 z_i - a_3 y_i + b \alpha_i, \quad i = 1,\ldots,N \tag{10.16}$$

$$0 = -a_1 z_i + a_3 x_i + b \beta_i, \quad i = 1,\ldots,N \tag{10.17}$$

$$0 = a_1 y_i - a_2 x_i + b \gamma_i, \quad i = 1,\ldots,N \tag{10.18}$$

Or in a more compact form

$$2\boldsymbol{r}_i = \|A\| \boldsymbol{v}_i, \quad b\boldsymbol{v}_i = \|A\| \boldsymbol{r}_i, \quad i = 1,\ldots,N \tag{10.19}$$

where the matrix $\|A\|$ has the form

$$\|A\| = \begin{Vmatrix} 0 & a_3 & -a_2 \\ -a_3 & 0 & a_1 \\ a_2 & -a_1 & 0 \end{Vmatrix}$$

Whence, taking into account the form of this matrix ("negative symmetry"), it follows, in particular, that the scalar product $(\boldsymbol{r}_i, \boldsymbol{v}_i) = 0$ for each $i = 1,\ldots,N$. That is, at any point $(\boldsymbol{r}, \boldsymbol{v}) \in B^*$, at which the minimum in (10.5) is reached under conditions (10.6) − (10.9), the velocity vector $\boldsymbol{v}_i$ of each particle of the system must be perpendicular to its radius-vector $\boldsymbol{r}_i$. From (10.19) further follows

$$2b\boldsymbol{r}_i = \|A\|^2 \boldsymbol{r}_i, \quad 2b\boldsymbol{v}_i = \|A\|^2 \boldsymbol{v}_i, \quad i = 1,\ldots,N \tag{10.20}$$

where the matrix

$$\|A\|^2 = \begin{Vmatrix} -a_2^2 - a_3^2 & a_1 a_2 & a_1 a_3 \\ a_1 a_2 & -a_1^2 - a_3^2 & a_2 a_3 \\ a_1 a_3 & a_2 a_3 & -a_1^2 - a_2^2 \end{Vmatrix}$$

By virtue of (10.20), at any point $(\boldsymbol{r}, \boldsymbol{v}) \in B^*$ of the minimum, the radius-vector of each particle $\boldsymbol{r}_i$ and its velocity vector $\boldsymbol{v}_i$, $i = 1,\ldots,N$ in three-dimensional space must be eigenvectors of the matrix $\|A\|^2$ with the same eigenvalue $2b$. Whence, taking into account also the symmetry of the matrix $\|A\|^2$, then follows

$$\|A\|^2 = 2b\|E\| \tag{10.21}$$

where $\|E\|$ is the identity matrix. From (10.21) we obtain the following equations for the coefficients $a_1, a_2, a_3, b$

$$a_1 a_2 = a_1 a_3 = a_2 a_3 = 0 \tag{10.22}$$

$$a_1^2 + a_2^2 = a_1^2 + a_3^2 = a_2^2 + a_3^2 = -2b \tag{10.23}$$

Due to (10.22), only one of the coefficients $a_1, a_2, a_3$ can be different from zero, whence it follows that the system of equations (10.22), (10.23) and the matrix equation (10.21) in the three-dimensional case has no solution, except for the trivial solution $a_1 = a_2 = a_3 = b = 0$. But with such a trivial solution, from (10.13) − (10.15) follows $x_i = y_i = z_i = 0$ for each $i = 1,\ldots,N$, which is inconsistent with condition (10.8).

Thus, at any point $(r, v) \in B^*$ of the minimum, all vectors $r_i$, $v_i$, $i = 1,\ldots, N$ must be in the same plane and, taking into account the conditions (10.6) $-$ (10.8) for the vector $L$, in the plane $z = 0$. Consequently, the coordinates of these vectors

$$z_i = 0, \quad \gamma_i = 0 \quad \text{for each} \quad i = 1,\ldots, N \tag{10.24}$$

In this case, the restrictions (10.6), (10.7) on the components $L_x, L_y$ of the moment of impulse are fulfilled automatically, that is coefficients $a_1 = a_2 = 0$. The initial problem of finding the minimum (10.5) under conditions (10.6) $-$ (10.9) in three-dimensional space becomes a problem in two-dimensional space on the plane $z = 0$.

The necessary minimum conditions (10.13) $-$ (10.18) for the coordinates of the vectors $r_i = (x_i, y_i)$, $v_i = (\propto_i, \beta_i)$, $i = 1,\ldots, N$ on this plane are then written in the form

$$2x_i = a_3 \beta_i, \quad 2y_i = -a_3 \propto_i, \quad i = 1,\ldots, N$$

$$b \propto_i = a_3 y_i, \quad b\beta_i = -a_3 x_i, \quad i = 1,\ldots, N$$

Whence follows the equality for the coefficients $2b = -a_3^2$, after which the necessary minimum conditions, taking into account (10.24), take the form

$$a_1 = a_2 = 0, \quad b = -a_3^2/2 \tag{10.25}$$

$$z_i = \gamma_i = 0 \quad \text{for each} \quad i = 1,\ldots, N \tag{10.26}$$

$$\propto_i = -(2/a_3)y_i, \quad \beta_i = (2/a_3)x_i, \quad i = 1,\ldots, N \tag{10.27}$$

Conditions (10.8), (10.9), taking into account these expressions, are written as

$$\sum_i m_i(x_i \beta_i - y_i \propto_i) = \frac{2}{a_3} \sum_i m_i(x_i^2 + y_i^2) = L \tag{10.28}$$

$$\frac{1}{2}\sum_i m_i(\propto_i^2 + \beta_i^2 + \gamma_i^2) = \frac{2}{a_3^2} \sum_i m_i(x_i^2 + y_i^2) = T \tag{10.29}$$

from which follows the equality for the coefficient $a_3 = L/T$.

From relations (10.24) $-$ (10.27), taking into account the obtained equalities for the coefficients $a_1$, $a_2$, $a_3$, $b$, it follows that the necessary conditions of minimum in problem (10.5) $-$ (10.9) at the point $(r, v) \in B^*$ have the form

$$a_1 = a_2 = 0, \quad a_3 = L/T, \quad b = -a_3^2/2 = -(L/T)^2/2 \tag{10.30}$$

$$z_i = \gamma_i = 0, \quad i = 1,\ldots, N \tag{10.31}$$

$$\propto_i = -(2T/L)y_i, \quad \beta_i = (2T/L)x_i, \quad i = 1,\ldots, N \tag{10.32}$$

where, due to (10.28), (10.29) the coordinates of the vectors $r_i = (x_i, y_i)$, $i = 1,\ldots, N$ must satisfy the relation

$$\sum_i m_i(x_i^2 + y_i^2) = \frac{L^2}{2T} \tag{10.33}$$

Whence, taking into account also the form of the objective function $g(r)$ in this problem, it follows that the minimum (10.5) on the set $B^*$ is achieved at any point $(r, v) \in B^*$ that satisfies the conditions (10.31) - (10.33), and this minimum is equal

$$\min g(r) = L^2/(2T) \tag{10.34}$$

For some points $(r, v) \in B^*$ that satisfy the necessary minimum conditions (10.31) − (10.33), the radius-vectors $r_i = (x_i, y_i)$ of different particles can coincide, which is impossible due to physical conditions in the initial phase space $\Omega$. Such points belong to the closure $B^*$, but do not belong to the original subset $B \subset \Omega$, over which the minimum is taken in problem (10.5) − (10.9). But as can be seen from these expressions, among the points that satisfy these necessary conditions, there are also points $(r, v) \in B \subset \Omega$ (that is, points with different vectors $r_i$) Therefore, in this problem the minimum (10.34) is also reached on the set $B \subset \Omega$. In other words, the minimum in problem (10.5) − (10.9), or in the original problem (10.2) − (10.4) on the phase space $\Omega$ is reached and equal to $L^2/2T = |L|^2/2T$. In addition, taking into account (10.24), (10.31), the solution of this problem can be achieved only on flat trajectories of the system. ∎

Proceeding from the previous Theorem 10.1, in the next section 11 it is shown (Theorem 11.1) that the minimum root-mean-square size (and the maximum compactness in this sense) of the system when it moves in three-dimensional space can be achieved only on flat trajectories.

**Theorem 10.2** Let $r_\lambda(t)$, $t \geq 0$ be a flat stationary and periodic trajectory from the family of trajectories (3.2) with characteristics $L \neq 0$, $T > 0$. Then the minimum in problem (10.5) − (10.9) (or in problem (10.2) − (10.4)) is reached at each point of this trajectory.

*Proof.* Due to the necessary minimum conditions (10.32) in any state $(r, v) \in B \subset \Omega$, in which the minimum is reached in problem (10.5) - (10.9), the velocity vector of each particle $v_i$ is perpendicular to its radius vector $r_i$. And the modulus of the velocity of each particle is equal to

$$|v_i| = (2T/L) |r_i| \quad \text{for each} \quad i = 1, \ldots, N \tag{10.35}$$

From the definition of the flat trajectory $r_\lambda(t)$, it follows that conditions (10.31) − (10.33) and (10.35) are satisfied at each point of this trajectory. In addition, taking into account relations (4.13), (4.14) obtained above in Section 4, relations for the characteristics of the system $g = L^2/(2T)$, $\omega = \sqrt{2\lambda} = 2T/L$ are satisfied on the trajectory $r_\lambda(t)$ (at each of its points), which corresponds to conditions (10.32) - (10.35). Thus, the minimum in (10.5) under conditions (10.6) − (10.9) is reached at each point of the flat stationary trajectory $r_\lambda(t)$. ∎

**Remark 1** on Theorems 10.1 and 10.2. By virtue of relations (4.15), for a given system with a set of masses of particles $m = (m_1, m_2, \ldots, m_N)$ on flat trajectories of the type $r_\lambda(t)$, not any values of the characteristics $T, L$ are possible, but only such that the equality $TL^2 = C^2(m)/2$ is satisfied (see Section 4 above). Theorem 10.1 is valid for any parameters $L \neq 0$, $T > 0$ of the problem (10.5) – (10.9). But if the parameters $TL^2 \neq C^2(m)/2$, then it is not known whether there are flat trajectories for the given system such that the minimum in this problem is reached at each point of the trajectory. Theorem 10.2 was obtained for trajectories $r_\lambda(t)$ for which the indicated equality holds automatically.

**Remark 2** on Theorem 10.2. The previous Theorem 10.2 was obtained for each trajectory $r_\lambda(t)$ generated by the state $r_\lambda$ of the global minimum in the optimization problem (2.5). But by virtue of Theorem 7.1, the main relations (4.13), (4.14) for the characteristics of the system are valid (up to the replacement of the function $C(m)$ by $C_k(m)$) for any trajectory $r_\lambda^{(k)}(t)$ generated by any local minimum $r_\lambda^{(k)}$ in problem (2.5) (see Sec. 7). Thus, theorem 10.2 is also valid for all flat stationary and periodic trajectories $r_\lambda^{(k)}(t)$ generated by local minima $r_\lambda^{(k)}$, $k = 1, \ldots, n$ in problem (2.5).

For the closed conservative system considered here, the angular momentum $L$ and the total energy $E = T + \Pi = T - f$, taking into account the known conservation laws ([1], [2] etc.), are constant. Let us consider the following problem of finding a minimum for the characteristics of a system in three-dimensional space: it is required to find

$$\min f(\boldsymbol{r}) \quad (10.36)$$

where the minimum is taken over all states of the phase space $(\boldsymbol{r}, \boldsymbol{v}) \in \Omega$ under conditions on the characteristics

$$g(\boldsymbol{r}) = g \quad (10.37)$$

$$\boldsymbol{L}(\boldsymbol{r}, \boldsymbol{v}) = \boldsymbol{L} \quad (10.38)$$

$$T(\boldsymbol{v}) - f(\boldsymbol{r}) = E \quad (10.39)$$

where $g$, $\boldsymbol{L} = (0, 0, L)$, $E$ are constants. This problem of finding the minimum in (10.36) under conditions (10.37) − (10.39) in three-dimensional space is a natural generalization of problem (2.4) considered above in Section 2 for the characteristics $f(\boldsymbol{r})$, $g(\boldsymbol{r})$ of the system on the plane.

**Theorem 10.3** *Let $r_\lambda(t)$, $t \geq 0$ be a flat stationary and periodic trajectory of a system from the family of trajectories (3.2) with characteristics $g$, $\boldsymbol{L} = (0, 0, L)$, $E$. Then the necessary conditions of minimum in (10.36) under conditions (10.37) − (10.39) are satisfied at each point of this trajectory.*

*Proof.* Let $P \subset \Omega$ be a subset of phase space over which the minimum in (10.36) is taken, that is, the set of all states $(\boldsymbol{r}, \boldsymbol{v}) \in \Omega$ in which conditions (10.37) − (10.39) are satisfied. It is easy to show that the minimum in (10.36) under conditions (10.37) − (10.39) is achieved on the set $P$. Denote by $P_c \subset P$ the subset of all points (states) $(\boldsymbol{r}, \boldsymbol{v}) \in P \subset \Omega$ of the phase space such that $f(\boldsymbol{r}) \leq c$, where $c$ is any sufficiently large constant such that $c > \inf f(\boldsymbol{r})$, where $\inf$ taken on the set $P \subset \Omega$. The function $f(\boldsymbol{r})$ is bounded on the set $P_c$, whence by virtue of condition (10.39) it follows that the function $T(\boldsymbol{v})$ is also bounded on this set. Whence, taking into account the condition (10.37) and the form of the functions $g(\boldsymbol{r})$, $T(\boldsymbol{v})$, it follows that the set $P_c \subset \Omega \subset E_{6N}$ is bounded. In addition, the functions $g(\boldsymbol{r})$, $\boldsymbol{L}(\boldsymbol{r}, \boldsymbol{v})$, $T(\boldsymbol{v})$ are continuous on the entire phase space $\Omega$, and the function $f(\boldsymbol{r})$ is continuous on the subset $P_c \subset P \subset \Omega$, that is the set $P_c$ is closed. Thus, the minimum (10.36) on a closed bounded subset $P_c \subset P$ of the phase space is achieved. Whence, taking into account the definition of this subset, it follows that the minimum (10.36) on the set $P$ is also achieved. The necessary conditions for this minimum at the point $(\boldsymbol{r}, \boldsymbol{v}) \in P$ are written as

$$-\text{grad } f(\boldsymbol{r}) = \sigma \text{ grad } g(\boldsymbol{r}) + \beta_1 \text{grad } L_x(\boldsymbol{r}, \boldsymbol{v}) + \beta_3 \text{grad } L_z(\boldsymbol{r}, \boldsymbol{v}) +$$

$$+\alpha \text{ grad } [T(\boldsymbol{v}) - f(\boldsymbol{r})] \quad (10.40)$$

where $\sigma$, $\beta_1$, $\beta_2$, $\beta_3$, $\alpha$ are the Lagrange multipliers. Equations (10.40) together with conditions (10.37) − (10.39) form a system of $(6N + 5)$ equations for $(6N + 5)$ unknowns, namely, for the point of minimum $\boldsymbol{r}_i = (x_i, y_i, z_i)$, $\boldsymbol{v}_i = (\alpha_i, \beta_i, \gamma_i)$, $i = 1, \ldots, N$; and coefficients $\sigma$, $\beta_1$, $\beta_2$, $\beta_3$, $\alpha$. Let us show that the necessary minimum conditions (10.40), (10.37) − (10.39) are satisfied at each point of the flat stationary trajectory $r_\lambda(t)$. By virtue of Theorem 3.2 (see Section 3 above), the following equality holds on this trajectory:

$$\text{grad } f[\boldsymbol{r}_\lambda(t)] = -\lambda \text{ grad } g[\boldsymbol{r}_\lambda(t)] \quad \text{for each} \quad t \geq 0$$

Where the gradient vectors, as in the previous equation (10.40), are considered in the phase space $\Omega$, (taking into account also the equalities $\partial g / \partial \boldsymbol{v}_i = 0$, $\partial f / \partial \boldsymbol{v}_i = 0$, $i = 1, \ldots, N$). Whence it follows that on the trajectory $r_\lambda(t)$, $t \geq 0$ the system of equations (10.40) has the form

$$(\lambda - \sigma - \alpha\lambda)\text{grad } g(\boldsymbol{r}) = \beta_1 \text{grad } L_x(\boldsymbol{r}, \boldsymbol{v}) + \beta_2 \text{grad } L_y(\boldsymbol{r}, \boldsymbol{v}) +$$

$$+\beta_3 \text{grad } L_z(\boldsymbol{r}, \boldsymbol{v}) + \alpha \text{ grad } T(\boldsymbol{v})$$

Let us determine the coefficient $\sigma$ from the equality $\sigma = \lambda$. Then the necessary minimum conditions (10.40) at each point of this trajectory are written as

$$\text{grad } g\,(r_\lambda) = a_1 \text{grad } L_x\,(r_\lambda, v) + a_2 \text{grad } L_y\,(r_\lambda, v) +$$
$$+ a_3 \text{grad } L_z\,(r_\lambda, v) + b\,\text{grad } T(v) \qquad (10.41)$$

where the coefficients

$$a_1 = -\beta_1/(\alpha\lambda),\ a_2 = -\beta_2/(\alpha\lambda),\ a_3 = -\beta_3/(\alpha\lambda),\ b = -1/\lambda \qquad (10.42)$$

Whence it follows that the conditions (10.41) and, thus, the necessary conditions for the minimum (10.40) in this problem on the trajectory $r_\lambda(t)$, $t \geq 0$ (at each of its points) coincide with the conditions (10.10) in the previous theorem for the problem (10.5) − (10.9) with up to the renaming in (10.42) of the Lagrange multipliers. Taking into account (10.40), (10.41), and repeating the proof of the previous theorem 10.1, we obtain, similarly to expressions (10.25) − (10.27), that the necessary minimum conditions in this problem on the trajectory $r_\lambda(t)$, $t \geq 0$ have the form

$$a_1 = a_2 = 0,\quad b = -a_3^2/2 \qquad (10.43)$$

$$z_i = \gamma_i = 0,\ i = 1, \ldots, N \qquad (10.44)$$

$$\alpha_i = -(2/a_3)y_i,\ \beta_i = (2/a_3)x_i,\ i = 1, \ldots, N \qquad (10.45)$$

Conditions (10.44) are satisfied on each flat stationary trajectory $r_\lambda(t)$, $t \geq 0$. In addition, from the definition of the trajectory $r_\lambda(t)$, $t \geq 0$ (see Section 3 above), it follows that conditions (10.45) are also satisfied on this trajectory (at each of its points), if the coefficient $a_3$ satisfies the relation

$$\omega = \sqrt{2\lambda} = 2/a_3 \qquad (10.46)$$

where $\omega = \sqrt{2\lambda}$ is the angular velocity of rotation of the system on the trajectory $r_\lambda(t)$. Then due to (10.43), (10.46) the equality in (10.42) for the Lagrange multiplier $b = -1/\lambda$ also holds. Note that, by virtue of (4.15) on any flat stationary trajectory $r_\lambda(t)$ the equality $f = 2T$ holds (see Section 4 above) and therefore the total energy of the system $E = T - f = -T$. Whence, taking into account also (4.11), (4.12), it follows that on the trajectory $r_\lambda(t)$ the relations

$$L = \sqrt{2\lambda}\,g,\ |E| = T = \lambda g, \qquad (10.47)$$

are satisfied. The condition (10.37) on the trajectory $r_\lambda(t)$, $t \geq 0$ (at each of its points) has the form

$$\sum_i m_i(x_i^2 + y_i^2) = const = g$$

The condition (10.38) at the angular momentum of the system on the trajectory $r_\lambda(t)$, $t \geq 0$ (at each of its points), taking into account (10.45), (10.46), is written as

$$L = \sum_i m_i(x_i\beta_i - y_i\alpha_i) = \sqrt{2\lambda}\sum_i m_i(x_i^2 + y_i^2) = \sqrt{2\lambda}\,g$$

The condition (10.39) on the total energy of the system on the trajectory $r_\lambda(t)$ due to (10.44) − (10.46), taking into account that the equality $E = -T$ is satisfied on this trajectory, is written as

$$|E| = T = \frac{1}{2}\sum_i m_i(\alpha_i^2 + \beta_i^2 + \gamma_i^2) = \lambda\sum_i m_i(x_i^2 + y_i^2) = \lambda g$$

which coincides with the conditions (10.47) for these characteristics on the trajectory $r_\lambda(t)$ for a given trajectory parameter $\lambda > 0$. Thus, this trajectory satisfies (at each of its points) the necessary minimum conditions in this problem. ∎

Note that Theorem 10.3 also remains valid (see the remark on Theorem 10.2 above) for all flat stationary and periodic trajectories $r_\lambda^{(k)}(t)$, generated by local minima $r_\lambda^{(k)}$ in problem (2.5).

Thus, the solution of the problems considered in this section for finding a conditional minimum for the characteristics of a system in three-dimensional space is achieved (for each $t \geq 0$) on the flat stationary and periodic trajectories $r_\lambda(t)$ and $r_\lambda^{(k)}(t)$ considered above.

### 11. Minimum possible size of the system can be achieved only on flat trajectories

Further, proceeding from Theorem 10.1 of the previous section, it is shown that the minimum possible root-mean-square size (and the maximum compactness in this sense) of the system when it moves in three-dimensional space can be achieved only on flat trajectories (Theorem 11.1). It is shown that such trajectories are, in particular, the flat stationary and periodic trajectories of the form $r_\lambda(t)$ and $\tilde{r}_\lambda(t, \omega_0)$ considered above, generated by the global minimum in the optimization problems considered in section 2.

Let $[r(t), v(t)] \in \Omega$, $t \geq 0$ be any stationary (in the sense of the definition of section 1) trajectory of the system in three-dimensional space and $b[r(t)]$ − root-mean-square size of the system on this trajectory at the moment of time $t \geq 0$, related to the function $g[r(t)]$ by the same relation as in the problem on the plane, i.e. $g[r(t)] = b^2[r(t)] \sum_j m_j$. Let $W(L, T)$ be the solution of the problem (considered above in section 10) of finding a conditional minimum in (10.2) − (10.4), that is

$$W(L, T) = \min g(r) \tag{11.1}$$

where the minimum is calculated over all states $(r, v) \in \Omega$ of the phase space under the conditions

$$L(r, v) = L, \quad T(v) = T \tag{11.2}$$

The following theorem shows that the minimum possible root-mean-square size (and the maximum compactness in this sense) of a system when it moves in three-dimensional space can be achieved only on flat trajectories.

**Theorem 11.1** *If at least one point of the stationary trajectory $[r(t), v(t)] \in \Omega$, $t \geq 0$ achieves the minimum possible root-mean-square size $b = b[r(t)]$ of the system for the current value of its cohesion $f = f[r(t)]$, then this trajectory can only be flat.*

*Proof.* In accordance with the known conservation laws ([1], [2] etc.), the angular momentum $L$ of the system and its total energy $E = T - f$ on the trajectory $[r(t), v(t)] \in \Omega$, $t \geq 0$ are some constants, that is

$$T[v(t)] - f[r(t)] \equiv E \quad \text{for each} \quad t \geq 0 \tag{11.3}$$

where $f[r(t)]$ is the cohesion and $T[v(t)]$ is the kinetic energy of the system at the time $t \geq 0$ on this trajectory, $E$ is a constant. From the definition of the problem of finding a conditional minimum (10.2) − (10.4) and relations (11.1), (11.2) it follows that on the trajectory $[r(t), v(t)] \in \Omega$ the function $g[r(t)] = g[r_1(t), r_2(t), ..., r_N(t)]$ satisfies the inequality

$$g[r(t)] \geq W\{L, T[v(t)]\} \quad \text{for each} \quad t \geq 0 \tag{11.4}$$

But according to the condition of this theorem (taking into account that $g[r(t)] = b^2[r(t)] \sum_j m_j$) on this trajectory the following equality is satisfied

$$g[r(t)] = W\{L, T[v(t)]\} \quad \text{for some} \quad t \geq 0$$

Where the kinetic energy $T[\boldsymbol{v}(t)]$ of the system at each current moment of time $t \geq 0$ due to identity (11.3) is uniquely determined by the value of its cohesion $f[\boldsymbol{r}(t)]$ at this moment. Then the proof follows from the previous Theorem 10.1 (see section 10 above). ∎

**Remark 1** on Theorem 11.1. The function $f = f[\boldsymbol{r}(t)]$ coincides with the potential energy of the system $\Pi = \Pi[\boldsymbol{r}(t)]$ up to a sign and an arbitrary constant, i.e. $f[\boldsymbol{r}(t)] = -\Pi[\boldsymbol{r}(t)] + C$. Thus, Theorem 11.1 can also be formulated in terms of the potential energy of the system in the following equivalent form: *If at least one point of the stationary trajectory $[\boldsymbol{r}(t), \boldsymbol{v}(t)] \in \Omega$, $t \geq 0$ achieves the minimum possible root-mean-square size $b = b[\boldsymbol{r}(t)]$ of the system for the current value of its potential energy $\Pi = \Pi[\boldsymbol{r}(t)]$, then this trajectory can only be flat.*

**Remark 2** on Theorem 11.1. By virtue of equality (11.3) Theorem 8.4 can also be formulated in terms of the kinetic energy of the system in the following equivalent form: *If at least one point of the stationary trajectory $[\boldsymbol{r}(t), \boldsymbol{v}(t)] \in \Omega$, $t \geq 0$ achieves the minimum possible root-mean-square size $b = b[\boldsymbol{r}(t)]$ of the system for the current value of its kinetic energy $T = T[\boldsymbol{v}(t)]$, then this trajectory can only be flat.*

Note that Theorem 11.1 remains valid if the stationarity condition for the trajectory $[\boldsymbol{r}(t), \boldsymbol{v}(t)] \in \Omega$, $t \geq 0$ is replaced by a weaker condition for its regularity in the following sense: A trajectory $[\boldsymbol{r}(t), \boldsymbol{v}(t)] \in \Omega$, $t \geq 0$ satisfying the system of differential equations (1.4) under some initial conditions (1.5) will be called regular, if there exists a constant $0 < C < \infty$ such that $|\boldsymbol{r}_j(t) - \boldsymbol{r}_i(t)| > C$ for each $t \geq 0$, $i \neq j$. Any stationary (in the sense indicated above in section 1) trajectory is obviously also regular, but the converse is generally not true.

Taking into account the definition of the function $W(\boldsymbol{L}, T)$ and relation (11.3), we will say that on the stationary trajectory $[\boldsymbol{r}(t), \boldsymbol{v}(t)] \in \Omega$ at the moment of time $t \geq 0$ the minimum possible root-mean-square size (and the maximum compactness in this sense) of the system is achieved, if at this moment the following equality is satisfied

$$g[\boldsymbol{r}(t)] = W\{\boldsymbol{L}, T[\boldsymbol{v}(t)]\}$$

where $\boldsymbol{L}$ is the angular momentum of the system on this trajectory and $T[\boldsymbol{v}(t)]$ is the kinetic energy of the system on this trajectory at the moment $t$.

It follows from the previous Theorems 10.1 and 10.2 that on each trajectory of the form $\boldsymbol{r}_\lambda(t)$ inequality (11.4) becomes the identity equality:

$$g[\boldsymbol{r}(t)] \equiv W\{\boldsymbol{L}, T[\boldsymbol{v}(t)]\} = L^2/(2T) \quad \text{for each} \quad t \geq 0 \quad (11.5)$$

where $L = |\boldsymbol{L}|$ is the value of angular momentum and $T \equiv T[\boldsymbol{v}(t)]$ is the kinetic energy of the system on this trajectory. Thus, on each trajectory of the form $\boldsymbol{r}_\lambda(t)$ the minimum possible root-mean-square size (maximum compactness) of the system is achieved (for each $t \geq 0$) when it moves in three-dimensional space.

We also consider flat trajectories of a more general form $\tilde{\boldsymbol{r}}_\lambda(t, \omega_0)$ generated by the global minimum in (2.5) (see Section 8 above). Taking into account formulas (10.34), (8.18), (8.21), (4.11), (4.16), after simple transformations, we find that on each trajectory of the form $\tilde{\boldsymbol{r}}_\lambda(t, \omega_0)$ the following relation is satisfied

$$\frac{g[\boldsymbol{r}(t)]}{W\{\boldsymbol{L}, T[\boldsymbol{v}(t)]\}} = G(\lambda, \omega_0, \varphi) \quad (11.6)$$

where $\varphi = \varphi(t)$ (see Section 8 above), and

$$G(\lambda, \omega_0, \varphi) = \frac{1 + e^2 - 2e\cos\varphi}{(1 - e\cos\varphi)^2}$$

where $e = 1 - (\omega_0^2/2\lambda)$. Whence it follows that the value (11.6) satisfies the following relations

$$G(\lambda, \omega_0, \varphi) \geq 1 \quad \text{for each } \varphi$$

$$G(\lambda, \omega_0, \varphi) = 1, \text{ if } \varphi = k\pi, \ k = 0,1,2,\ldots$$

Thus, on each flat trajectory of the form $\tilde{r}_\lambda(t, \omega_0)$ the minimum possible root-mean-square size of the system during its motion in three-dimensional space is achieved at the points of the trajectory with the angular coordinate $\varphi = k\pi, \ k = 0,1,2,\ldots$. In the particular case when the parameters of the trajectory are $\omega_0 = \sqrt{2\lambda}$ (which corresponds to trajectories of the form $r_\lambda(t)$), then $G(\lambda, \omega_0, \varphi) \equiv 1$ and the minimum possible size of the system is achieved at each point of the trajectory.

Thus, due to Theorems 10.1, 11.1, the minimum possible root-mean-square size (and the maximum compactness in this sense) of the system when it moves in three-dimensional space can be achieved only on flat trajectories. And such trajectories are, in particular, flat stationary and periodic trajectories of the form $r_\lambda(t)$ and $\tilde{r}_\lambda(t, \omega_0)$ generated by the global minimum in the optimization problem (2.5) (see also Theorems 4.3 and 8.4 above).

# REFERENCES

[1]  L. D. Landau, E. M. Lifshits, Mechanics, Moscow: Science, 1965, 204 pp.
[2]  F. R. Gantmacher, Lectures on analytical mechanics, Moscow: Science, 2001, 264 pp.
[3]  C. L. Siegel, Lectures on celestial mechanics, Berlin, New York: Springer, 1959, 290 pp.
[4]  W.I. Zangwill, Nonlinear programming, Prentice – Hall, Inc. Englewood, 1969, 312 pp.
[5]  N. N. Moiseev, Numerical methods in theory of optimal systems, Moscow: Science, 1971, 424 pp.
[6]  N. N. Moiseev, Theory of optimal systems, Moscow: Science, 1975, 528 pp.
[7]  M. Hampton, R. Moeckel, Finiteness of relative equilibria of the four-body problem, *Inventions mathematics*, **163**, 2006, 289 – 312.
[8]  J. Montaldie, Existence of symmetric central configurations. *Celestial Mechanics and Dynamical Astronomy*, **122**, 2015, 405 – 418.
[9]  E. Roy, Orbital Motion, CRC Press, 2020, 544 pp.
[10] V.I. Arnold, Small denominators and problems of movement stability in classical and celestial mechanics, Successes of mathematical sciences, **18** (6), 1963, 91–192.
[11] G. D. Birkhoff, Dynamical systems, New York, American Mathematical Society, 1927, 316 pp.

Igor Pavlov, Department of Mathematics, Bauman Moscow State Technical University, ul. Baumanskaya 2-ya, 5, Moscow, 105005, Russian Federation.

E-mail address: ipavlov@bmstu.ru